\documentclass[11pt,a4paper]{article}

 \usepackage{amssymb,amsmath,amsfonts,amssymb}
 \usepackage{graphics,graphicx,color}
 -\marginparwidth \oddsidemargin -4mm \evensidemargin
 \oddsidemargin%
 \advance\hoffset by 5mm
 \usepackage{amsmath}    
 \advance\hoffset by 5mm

\newcommand\grad{{\bf \nabla}}

\newcommand\la{{\lambda}}

\newcommand\n{{\bf n}}
\newcommand\zhat{{\bf z}}
\newcommand\x{{\bf x}}
\newcommand\e{{\bf e}}
\newcommand\R{{\bf R}}
\newcommand\y{{\bf y}}

\newcommand\m{{\bf m}}
\newcommand\I{{\bf I}}
\newcommand {\Ical}{{\cal{{I}_{LG}}}}

\newcommand\p{{\bf p}}

\newcommand\eps{{\epsilon}}
\newcommand\Q{{\bf Q}}
\newcommand{\Acal}{{\cal A}}

\newcommand{\Rr}{{\mathbb R}}

\newcommand{\sgn}{\,\text{sgn}\,}

 \title{Uniaxiality in the Landau-de Gennes theory of nematic liquid
crystals}
 \date{\today}
\author{Apala Majumdar\footnote{ Mathematical Institute, University of Oxford, 24-29 St. Giles', OX1 3LB, U.K.}}

 \newtheorem{lemma} {Lemma}
 \newtheorem{proposition}{Proposition}

\begin{document}
 \maketitle
\begin{abstract}
We study uniaxial energy minimizers within the Landau-de Gennes
theory for nematic liquid crystals, subject to dirichlet boundary
conditions. Topological defects in such minimizers correspond to
the zeros of the corresponding equilibrium field. We consider
two-dimensional and three-dimensional domains separately and study
the correspondence between Landau-de Gennes theory and
Ginzburg-Landau theory for superconductors. We obtain results for
the  location and dimensionality of the defect set, the minimizer
profile near the defect set and study the qualitative properties
of uniaxial energy minimizers away from the defect set, in the
physically relevant case of vanishing elastic constant. In the
three-dimensional case, we establish the
$C^{1,\alpha}$-convergence of uniaxial minimizers to a limiting
harmonic map, away from the defect set, for some $0<\alpha<1$.
Some generalizations for biaxial minimizers are also discussed.
This work is motivated by the study of defects in liquid
crystalline systems and their applications.
\end{abstract}

\section{Introduction}
\label{sec:intro}

Nematic liquid crystals are examples of \emph{mesophases} whose physical
properties are intermediate between those of a typical liquid
and a crystalline solid \cite{dg}. The constituent rod-like
molecules have no translational order but exhibit a certain degree
of long-range orientational ordering. Consequently, liquid
crystals are anisotropic media and this makes them suitable for a
wide range of physical applications and the subject of very
interesting mathematical modelling \cite{lin}.

The Landau-de Gennes theory is a general continuum theory for
nematic liquid crystals \cite{dg,virga}. It describes the state of
a nematic liquid crystal by a symmetric, traceless $3\times 3$
matrix - the $\Q$-tensor order parameter, that is defined in terms
of anisotropic macroscopic quantities, such as the magnetic
susceptibility and the dielectric anisotropy. Nematic liquid
crystals are said to be in the (a) \emph{biaxial} phase when $\Q$
has three distinct eigenvalues, (b) \emph{uniaxial} phase when
$\Q$ has a pair of equal non-zero eigenvalues and (c)
\emph{isotropic} phase when $\Q$ has three equal eigenvalues or
equivalently when $\Q=0$. For a general biaxial phase, $\Q$ can be
written in the form \cite{amaz, newtonmottram}
\begin{equation}
\label{eq:1} \Q =s\left(\n\otimes \n-\frac{1}{3}\I\right) +
r\left(\m\otimes \m-\frac{1}{3}\I\right) \quad s,r \in \Rr;~ \n,\m
\in S^2,
\end{equation} where $s,r$ are scalar order parameters, $\n,\m$
are eigenvectors of $\Q$ and $\I$ is the $3\times 3$ identity
matrix. In the uniaxial phase, $\Q$ takes the simpler form of
\begin{equation}
\label{eq:2} \Q =s\left(\n\otimes \n-\frac{1}{3}\I\right)  \quad s
\in \Rr;~ \n \in S^2
\end{equation} where $s$ is a scalar order parameter that measures
the degree of orientational ordering about the distinguished
eigenvector $\n$.

The Landau-de Gennes energy functional, $\Ical$, is a nonlinear
integral functional of $\Q$ and its spatial derivatives. In the
absence of any surface energies or external fields, $\Ical$ is
given by \cite{dg,newtonmottram}
\begin{equation}
\label{eq:3} \Ical[\Q] = \int_{\Omega} \frac{1}{2}|\grad \Q|^2 +
\frac{f_B(\Q)}{L}~dV
\end{equation}
where $\Omega$ is the domain, $f_B(\Q)$ is the \emph{bulk} energy
density that dictates the preferred phase - isotropic, uniaxial or
biaxial, $L$ is a positive material-dependent elastic constant and
$|\grad \Q|^2$ is an \emph{elastic} energy density that penalizes
spatial inhomogeneities. The equilibrium, physically
observable configurations correspond to either global or local
Landau-de Gennes energy minimizers, subject to the imposed
boundary conditions.

In this paper, we study uniaxial global minimizers (of the form
(\ref{eq:2})) of the Landau-de Gennes energy functional. There is
the very important underlying question - do uniaxial global
minimizers actually exist \cite{amaz}? This is an open question
but there is substantial numerical and experimental evidence to
show that global energy minimizers are \emph{largely uniaxial}
almost everywhere, in the sense that they have a small degree of
biaxiality. Therefore, a rigorous study of uniaxial global
minimizers is the first step in the mathematical analysis of
arbitrary global minimizers and the interplay between biaxiality
and uniaxiality. Secondly, stable uniaxial configurations do
exist, at least for certain temperature regimes and certain
physically realistic choices of the material-dependent constants
\cite{mg,rossovirga,sonnet}. Some of our results extend to stable
uniaxial configurations with appropriately \textit{bounded} energy
and can then be used to understand these configuration structures and
the nature of their singularities. We also point out that there
are two widely-used continuum theories for purely uniaxial liquid
crystal phases - the \textit{Oseen-Frank} theory and the
\textit{Ericksen} theory, both of which have received considerable
attention amongst the mathematical analysts and the numerical
modellers \cite{bcl,partialcrystal,linpoon, ericksen}. In fact,
uniaxiality is one of the most frequently used assumptions in the
theoretical study of liquid crystalline systems, even in the
context of applications.

We assume that uniaxial global minimizers exist for each $L>0$ and
then establish various properties of such uniaxial global
minimizers throughout the paper. The paper is organized as
follows. In Section~\ref{sec:prelim}, we introduce some basic
notation and terminology. In Section~\ref{sec:2D}, we study
Landau-de Gennes minimizers on two-dimensional (2D) domains and
establish a $1-1$ correspondence between Landau-de Gennes theory
and Ginzburg-Landau theory. In Section~\ref{sec:results1}, we
recall useful results from \cite{amaz} that are crucial for the
development of this paper. We study uniaxial Landau-de Gennes
minimizers on three-dimensional (3D) domains in the
low-temperature regime. There are important differences between
the $2D$ and $3D$ cases and the standard Ginzburg-Landau
techniques do not extend to the $3D$ case. We derive the governing
equations for uniaxial global minimizers and obtain qualitative
information about topological defects. The scalar order parameter
$`s'$ necessarily vanishes at the defect locations.  The defect
locations are prescribed in terms of the singular set of a
\emph{limiting harmonic map} and using asymptotic methods, we show
that the leading eigenvector $\n$ (see (\ref{eq:2})) necessarily
has a radial-hedgehog type of profile in the immediate
neighbourhood of each isolated point defect. Our result is
analogous to a powerful result on singularity profiles in
\cite{bcl}, where the authors work within the Oseen-Frank theory
for uniaxial liquid crystals with constant order parameter $s$.
In Section~\ref{sec:results2}, we study the qualitative properties
of uniaxial global minimizers away from the defect set, in the
limit $L\to 0^+$. The elastic constant $L$ to typically several
orders of magnitude smaller than the other material-dependent
constants and hence the $L\to 0$ limit is physically realistic
\cite{elastic}. We adapt the small energy regularity theorem of
\cite{chen} to the Landau-de Gennes framework and prove the
$C^{1,\alpha}$-convergence of uniaxial global minimizers to a
limiting harmonic map, away from the defect set, as $L\to 0$. This
convergence result encodes quantitative information about the
corresponding scalar order parameter. In Section~\ref{sec:dis}, we
discuss various generalizations of our results to uniaxial
solutions with bounded energy and to the completely general
biaxial case. The uniaxial case in $3D$ can be viewed as a
generalized Ginzburg-Landau theory from $\Rr^3$ to $\Rr^3$
although there are important technical differences. However, the
biaxial case presents a whole host of new mathematical
difficulties; there are five degrees of freedom in the biaxial
case and the additional degrees of freedom give us more
possibilities, particularly in the context of defects. The methods
in this paper contribute to the development of a generalized
Ginzburg-Landau theory from $\Rr^3$ to higher dimensions ($\Rr^5$
in this case), for non-standard non-convex multi-well bulk
potentials.

\section{Preliminaries}
\label{sec:prelim}

Let $\bar{S}\subset \mathbb{M}^{3\times 3}$ denote the space of
symmetric, traceless $3\times 3$ matrices  i.e.
\begin{displaymath}
\label{eq:4}
 \bar{S}\stackrel{def}{=} \left\{\Q \in \mathbb{M}^{3\times 3};
\Q_{ij}=\Q_{ji},~\Q_{ii} = 0 \right\}
\end{displaymath}
where we have used the Einstein summation convention; the Einstein
convention will be used in the rest of the paper. The
corresponding matrix norm is defined to be
\begin{displaymath}
\label{eq:5}
 \left| \Q \right|\stackrel{def}{=}\sqrt{\textrm{tr}\Q^2} =\sqrt{ \Q_{ij}
\Q_{ij}} \quad i,j=1\ldots 3.
\end{displaymath}

We take our domain $\Omega$ to be either a two-dimensional or
three-dimensional bounded, connected and simply-connected set with
smooth boundary, $\partial\Omega$. We work with the simplest form
of the bulk energy density, $f_B$, in (\ref{eq:3}) that allows for
a first-order nematic-isotropic phase transition
\cite{newtonmottram}. We focus on the low-temperature regime; the
function $f_B$ is bounded from below and can be written as
\begin{equation}
\label{eq:6}f_B(\Q) = -\frac{a^2}{2}\textrm{tr}\left(\Q^2\right) -
\frac{b^2}{3}\textrm{tr}\left(\Q^3\right) +
\frac{c^2}{4}\left(\textrm{tr}(\Q^2)\right)^2 +C(a^2,b^2,c^2)
\end{equation} where $a^2,b^2,c^2\in \Rr^+$ are material-dependent
and temperature-dependent positive constants and $C(a^2,b^2,c^2)$
is a positive constant that ensures $f_B(\Q)\geq 0$ for all
$\Q$-tensors. We note that $C(a^2,b^2,c^2)$ plays no role in
energy minimization, in either spatially homogeneous or
inhomogeneous cases. The negative coefficient of
$\textrm{tr}\left(\Q^2\right)$ incorporates the fact that we are
working below the nematic-isotropic transition temperature where
$f_B$ attains its minimum on the set of uniaxial $\Q$-tensors
given by \cite{bm}
\begin{eqnarray}
\label{eq:7} && \Q_{min} =
 \left\{\Q\in \bar{S}, \Q= s_+ \left( \n\otimes \n -
\frac{1}{3}\I \right)~
\right\}\end{eqnarray} with $\n\in\mathbb{S}^2$ and
\begin{equation}  s_+=\frac{b^2+\sqrt{b^4+24a^2c^2}}{4c^2}.
\label{eq:8}
\end{equation}

We study uniaxial global minimizers of $\Ical$, with $f_B$ as in
(\ref{eq:6}), with \emph{strong anchoring conditions} or dirichlet
boundary conditions. The prescribed boundary condition $\Q_b$ is
given by
\begin{equation}
\label{eq:9} \Q_b = s_+\left(\n_b \otimes
\n_b - \frac{1}{3}\I\right)
\end{equation} where $\n_b\in W^{1,2}(\Omega;M)$ ($M=S^1$ in $2D$ and $M=S^2$ in $3D$)
is a unit-vector field with non-zero topological degree $d$, when viewed as a map from $\partial\Omega$ to $M$.
Clearly, $\Q_b \in \Q_{\min}$ where $\Q_{\min}$ has been defined
in (\ref{eq:7}). We define our admissible space to be
\begin{eqnarray}
&& \Acal_{\Q} = \left\{\Q\in W^{1,2}\left(\Omega;\bar{S}\right);
\textrm{$\Q=\Q_b$ on $\partial\Omega$, \textrm{ with $\Q_b$ as in
}(\ref{eq:9})}\right\}\label{eq:10},
\end{eqnarray} where $W^{1,2}\left(\Omega;\bar{S}\right)$ is the
Sobolev space of square-integrable $\Q$-tensors with
square-integrable first derivatives \cite{evans}. The existence of
global energy minimizers for $\Ical$ in the admissible space
$\Acal_\Q$ follows readily from the direct methods in the calculus
of variations \cite{bm,amaz}. For completeness, we recall that the
$W^{1,2}$-norm is given by $\| \Q \|_{W^{1,2}(\Omega)} =\left(
\int_{\Omega} |\Q|^2 + |\grad \Q|^2~dx\right)^{1/2}.$ In addition
to the $W^{1,2}$-norm, we also use the $L^{\infty}$-norm in this
paper, defined to be $\|\Q\|_{L^{\infty}(\Omega)} = \textrm{ess
sup}_{\x\in\Omega}|\Q(\x)|$ .

Finally, we introduce the concept of a  \textit{``limiting
uniaxial harmonic map''} $\Q^{0}:\Omega \to \Q_{min}$; $\Q^{0} $
is defined to be
\begin{equation}
\label{eq:Q0} \Q^{0} = s_+\left(\n_0 \otimes
\n_0 - \frac{1}{3}\I\right)\end{equation}  where $\n_0$ is a minimizer of the Dirichlet energy
\begin{equation}
\label{eq:dirichlet}
I_{OF}[\n] = \int_{\Omega} \left|\grad \n \right|^2~dV
\end{equation}
in the admissible space
\begin{equation}
\label{eq:dirichlet2} \Acal_\n = \left\{\n\in
W^{1,2}\left(\Omega;M\right);~ \n = \n_b~on~\partial \Omega
\right\}
\end{equation} where $M=S^1$ in $2D$ and $M=S^2$ in $3D$.
 The
terminology \emph{limiting harmonic map} stems from the fact that
$\n_0$ is a harmonic unit-vector field \cite {schoen} and it can
be shown that $\Q^0$ is a global minimizer of $\Ical$ in the
restricted class $\Acal_\Q \cap \left\{\Q_{\min}\right\}$
\cite{bcl},\cite{lin}. We use the limiting harmonic map $\Q^0$ to
study the inter-relationship between the Landau-de Gennes theory
and the Oseen-Frank theory for nematic liquid crystals. The
Oseen-Frank theory is the simplest continuum theory for nematic
liquid crystals, restricted to uniaxial phases with constant
scalar order parameter \cite{dg}. Working within the one-constant
approximation, the Oseen-Frank energy reduces to the Dirichlet
energy in (\ref{eq:dirichlet}) and $\n_0$, and hence $\Q^0$, is a
global Oseen-Frank energy minimizer in the admissible space
$\Acal_\n$.

\section{The $2D$ case}
\label{sec:2D}

Let $\Omega\subset \Rr^2$ be a bounded, connected and
simply-connected domain with smooth boundary. Let $\bar{S}_2$ denote the space of
symmetric, traceless $2\times 2$ matrices. Then $\Q\in \bar{S}_2$ can be written as
\begin{equation}
\label{2D1} \Q = \la\left(\n\otimes\n - \m\otimes \m\right)
\end{equation}
where $\la\in\Rr$ and $\n, \m$ are the two orthonormal
eigenvectors of $\Q$. We note that there are only two degrees of
freedom in the representation (\ref{2D1}) and hence we can think
of $\Q$ as being a map $\Q:\Omega \to \Rr^2$. Using the identity,
$\delta_{ij} = \n_i\n_j + \m_i\m_j$, we can re-write (\ref{2D1})
as
\begin{equation}
\label{2D2} \Q = 2\la\left(\n\otimes\n -
\frac{1}{2}\mathbf{I}\right)
\end{equation}
where $\mathbf{I}$ is the $2\times 2$ identity matrix. Thus, all
admissible $\Q$-tensors in two dimensions necessarily have a
uniaxial structure as in (\ref{eq:2}) \begin{footnote}{We can also think of $\Q\in\bar{S}_2$ as being a symmetric, traceless $3\times3$ matrix: $\Q = (\la + \frac{1}{6})\n\otimes \n + (\frac{1}{6}-\la)\m\otimes \m - \frac{1}{3}\zhat\otimes\zhat \in \bar{S}$ where $\zhat$ is the unit-vector in the $z$-direction and $\bar{S}$ is the space of symmetric, traceless $3\times 3$ matrices.}\end{footnote}.

Straightforward calculations show that
\begin{eqnarray}
\label{2D3} && \left| \Q \right|^2 = 2\la^2 \nonumber\\ &&
\textrm{tr}\Q^3 = \Q_{ij}\Q_{jp}\Q_{pi} = 0 \quad i,j,p=1,2.
\end{eqnarray} Then the Landau-de Gennes energy functional in
(\ref{eq:3}) simplifies to
\begin{equation}
\label{2D4} \Ical[\Q] = \int_{\Omega}\frac{1}{2}\left|\grad
\Q\right|^2 + \frac{1}{L}\left\{-\frac{a^2}{2}\textrm{tr}\Q^2 +
\frac{c^2}{4}(\textrm{tr}\Q^2)^2\right\}~dV
\end{equation}
for two-dimensional domains. The corresponding Euler-Lagrange
equations are -
\begin{equation}
\label{2D5} \Q_{ij,kk} = \frac{1}{L}\left(-a^2 + c^2
|\Q|^2\right)\Q_{ij} \quad i,j=1,2
\end{equation}
and using the scaling $\widetilde{\Q} = \sqrt{\frac{c^2}{a^2}}\Q$,
we obtain the following system of partial differential equations -
\begin{eqnarray}
\label{2D6} && \widetilde{\Q}_{ij,kk} =
\frac{a^2}{L}\left(\left|\widetilde{\Q}\right|^2 -
1\right)\widetilde{\Q}_{ij} \quad i,j,k=1,2 \nonumber\\
&& \widetilde{\Q} = 2\left(\n_b\otimes \n_b
-\frac{1}{2}\mathbf{I}\right)~on~\partial\Omega.
\end{eqnarray} This is identical to the Ginzburg-Landau equations
for superconductors in two dimensions \cite{shafrir2} and we are
interested in the asymptotic properties of global energy
minimizers either in the limit $a^2\to \infty$ or $L\to 0^+$.

Let $\Q^L$ be a global minimizer of $\Ical$ in (\ref{2D4}), in the
admissible space \newline $\Acal_\Q = \left\{\Q\in
W^{1,2}\left(\Omega;\bar{S}_2\right); \Q = s_+\left(\n_b\otimes
\n_b - \frac{1}{2}\mathbf{I}\right)~on~\partial\Omega \right\}$ for a fixed $L>0$.
Then $\Q^L$ is necessarily of the form
\begin{equation}
\label{2D7} \Q^L(\x) = s^L(\x)\left(\n^L(\x)\otimes \n^L(\x) -
\frac{1}{2}\mathbf{I}\right)
\end{equation}
for some scalar function $s^L:\bar{\Omega}\to \Rr$ and $\n^L\in
W^{1,2}\left(\Omega;S^1\right)$. Let $\Theta_L = \left\{
\x\in\Omega; s^L(\x) = 0\right\}$ denote the isotropic set of
$\Q^L$. We have a topologically non-trivial boundary condition in
(\ref{2D6}), since $\n_b$ has non-zero topological degree when
viewed as a map from $\partial\Omega$ to $S^1$. Hence, the
unit-vector field $\n^L$ necessarily has interior discontinuities
and let $S_\n$ denote the defect set of $\n^L$. Then
\begin{equation}
\label{2D8} S_\n \subset \Theta_L
\end{equation}
and in what follows, we use existing results in the mathematical
literature for Ginzburg-Landau theory in two dimensions, to make
predictions about the structure and location of the isotropic set
and the far-field properties of global energy minimizers.

\noindent \textbf{Dimension of $\Theta_L$ \cite{bauman,elliott}:}
The isotropic set $\Theta_L$ consists of $|d|$ isolated points,
$\left\{a_1,\ldots a_{|d|}\right\}$, where $d$ is the topological
degree of the boundary condition $\Q_b$ in (\ref{eq:9}).

\noindent \textbf{Defect locations \cite{bbh2,shafrir2}:} The
configuration $\left(a_1,\ldots a_{|d|}\right)$ minimizes the
renormalized energy $W$ over $(b_1,\ldots
b_{|d|})\in\Omega^{|d|}$, which is defined by
\begin{equation}
\label{2D9} W(b_1,\ldots b_{|d|}) = -2\pi\sum_{i\neq j} \log|b_i -
b_j| - 2\pi\sum_{i,j}R(b_i,b_j)
\end{equation}
where $R(\x,\y) = \Psi(\x,\y) - \log|\x-\y|, \x,\y\in \Rr^2$ and
$\Psi(\x,\y)$ is given by the solution of an explicit
boundary-value problem.

\noindent \textbf{Far-field behaviour \cite{bbh2,shafrir2}:} Let
$\left\{ \Q^{L_k}\right\}$ denote a sequence of global energy
minimizers for (\ref{2D4}), where $L_k\to 0^+$ as $k\to\infty$.
Then (up to a subsequence),
$$ \Q^{L_k} \to \Q^* ~in~ C^{1,\alpha}\left(\bar{\Omega}\setminus
\Theta_L \right), \quad \forall\alpha<1 ~and~in~W^{1,p}(\Omega),
~\forall p\in\left[1,2\right)$$ for some $\Q^*\in \cap_{1\leq p
<2} W^{1,p}\left(\Omega;S^1\right)$. The limit $\Q^*$ is the
canonical harmonic map associated with $a_1,\ldots a_{|d|}$ and
the degrees $\sgn d,\ldots \sgn d$.

\noindent The interested reader is referred to \cite{bbh2,
shafrir2} for the proofs.

\section{Uniaxial minimizers and their defect sets in $3D$}
\label{sec:results1}

Let $\Omega\subset \Rr^3$ be a bounded, connected and
simply-connected domain with smooth boundary. An arbitrary
$\Q$-tensor field, $\Q:\Omega\to\bar{S}$ can be written as
$$\Q = \sum_{i=1}^{3}\la_i \e_i\otimes\e_i \quad
\sum_{i}\la_i=0$$ where $\e_i$ are the orthonormal eigenvectors,
$\la_i$ are the corresponding eigenvalues and $\textrm{tr}\Q^3\neq
0$ in general.

We study uniaxial global minimizers of the Landau-de Gennes energy
functional, $\Ical$ in (\ref{eq:3}), in the admissible space
$\Acal_\Q = \left\{\Q\in W^{1,2}\left(\Omega;\bar{S}\right); \Q =
\Q_b~on~\partial\Omega \right\}$. The corresponding Euler-Lagrange
equations are
\begin{equation}
L\Delta
\Q_{ij}=-a^2\Q_{ij}-b^2\left(\Q_{ik}\Q_{kj}-\frac{\delta_{ij}}{3}\textrm{tr}(\Q^2)\right)
+c^2\Q_{ij}\textrm{tr}(\Q^2)\,~ ~i,j=1,2,3 \label{eq:11},
\end{equation} where the term $b^2\frac{\delta_{ij}}{3}\textrm{tr}(\Q^2)$ is a Lagrange multiplier
associated with the tracelessness constraint. It follows from
standard arguments in elliptic regularity that a global minimizer $\Q^{*}$ is
actually a classical solution of (\ref{eq:11}) and $\Q^{*}$ is
smooth and real analytic on $\Omega$, up to the boundary
\cite{amaz}.

We assume that a uniaxial global minimizer exists for each $L>0$.
For $\Q$ uniaxial (of the form $\Q=s(\n\otimes\n -
\frac{1}{3}\mathbf{I})$ where $s:\Omega\to \Rr$ and $\n\in
W^{1,2}(\Omega;S^2)$, see (\ref{eq:2})), a direct calculation
shows that
$$\left(\Q_{ik}\Q_{kj}-\frac{\delta_{ij}}{3}\textrm{tr}(\Q^2)\right)
= \frac{s}{3}\Q_{ij}$$ and hence, the Euler-Lagrange equations
(\ref{eq:11}) simplify to
\begin{equation}
\label{eq:decoupled} L\Q_{ij,kk} = \frac{1}{3}\left(2c^2 s^2 -
b^2s - 3a^2\right)\Q_{ij}, \quad i,j=1\ldots 3.
\end{equation} Let $\Q^L$ denote a uniaxial global Landau-de
Gennes minimizer for a fixed $L>0$. Then $\Q^L$ is a classical
solution of (\ref{eq:decoupled}) and we are interested in the
qualitative properties of $\Q^L$ in the limit $L\to 0$.

We briefly comment on limiting harmonic maps in a $3D$ setting:
$\Q^0 = s_+\left(\n_0\otimes \n_0 - \frac{1}{3}\mathbf{I}\right)$
where $\n_0$ is an energy minimizing harmonic map in the
admissible space $\Acal_\n = \left\{\n\in W^{1,2}(\Omega;S^2);
\n=\n_b~on~\partial\Omega\right\}$. Let $S_0$ denote the singular
set of $\n_0$ (and hence, of $\Q^0$). Then $S_0$ consists of
precisely $|d|$ isolated point singularities \cite{bcl,schoen}.

 We, next, quote important results from
\cite{amaz, bm} that are crucial for the analysis in this paper.

\noindent \textbf{Maximum principle \cite{bm}:} Let $\Q$ be an
arbitrary solution (not necessarily uniaxial) of the
Euler-Lagrange equations (\ref{eq:11}) in the space $\Acal_\Q$.
Then
\begin{equation}
\label{eq:max} \|\Q\|_{L^{\infty}(\Omega)} \leq
\sqrt{\frac{2}{3}}s_+
\end{equation}
where $s_+$ has been defined in (\ref{eq:8}).

\noindent \textbf{Strong convergence to $\Q^0$ \cite{amaz}:} Let
$\Omega\subset \Rr^3$ be a bounded, connected and simply-connected domain with
smooth boundary. Let $\left\{\Q^{L_k}\right\}$ be a sequence of uniaxial global
minimizers of $\Ical$ in the admissible space $\Acal_\Q$ ( $\Ical$
and $\Acal_\Q$ have been defined in (\ref{eq:3}) and (\ref{eq:10})
respectively) where $L_k \to 0$ as $k\to\infty$. Then $\Q^{L_k}
\to \Q^0$ strongly in $W^{1,2}\left(\Omega;\bar{S}\right)$ (upto a subsequence), where
$\Q^0$ has been defined in (\ref{eq:Q0}).

\noindent \textbf{Interior and boundary monotonicity lemmas
\cite{amaz}:} Let $\Q$ be an arbitrary solution of the
Euler-Lagrange equations (\ref{eq:11}). Define the normalized
energy on balls $B(\x,r)\subset \Omega = \left\{ \y \in \Omega:
|\x - \y|\leq r\right\}$:
$$ \mathcal{F}(\Q,\x,r) = \frac{1}{r}\int_{B(\x,r)}\frac{1}{2}|\grad \Q|^2 +
\frac{f_B(\Q)}{L}~dV.$$ Then we have the following interior
monotonicity lemma:
\begin{equation}
\label{eq:intmon} \mathcal{F}(\Q,\x,r) \leq \mathcal{F}(\Q,\x,R)
\quad \forall \x \in \Omega;~ r\leq R ~and ~ B(\x,R)\subset
\Omega.
\end{equation}
Similarly, for $\x_0\in\partial\Omega$, we define the region
$\Omega_r = \Omega \cap B(\x_0,r)$ with $r>0$, and the
corresponding normalized energy to be
$$ \mathcal{E}(\Q,\x_0,r) = \frac{1}{r}\int_{\Omega_r}\frac{1}{2}|\grad \Q|^2 +
\frac{f_B(\Q)}{L}~dV.$$ Then there exists $r_0 > 0$ so that
\begin{equation}
\label{eq:bnd} \frac{d}{dr}\mathcal{E} \geq -
C\left(a^2,b^2,c^2,\Q_b,r_0,\Omega\right) \quad 0<r<r_0
\end{equation}
where the positive constant $C$ is independent of $L$.

The proofs of (\ref{eq:intmon}) and (\ref{eq:bnd}) follow a
standard pattern using the Pohozaev identity; complete details can
be found in \cite{amaz}. An immediate consequence of the strong
convergence and the monotonicity lemmas is the following:

\noindent \textbf{Convergence of bulk energy density away from
$S_0$ \cite{amaz}:} Let $\left\{\Q^{L_k}\right\}$ be a sequence of
global Landau-de Gennes energy minimizers in the admissible space
$\Acal_\Q$, where $L_k\to 0$ as $k \to \infty$. Assume that we
have a sequence $\Q^{L_k}$ with $L_k\to 0$ as $k\to \infty$, such
that $\Q^{L_k}\to \Q^0$ in $W^{1,2}(\Omega,\bar{S})$, as
$k\to\infty$, where $\Q^0$ has been defined in (\ref{eq:Q0}).

For any compact set $K\subset \bar{\Omega}$ such that $K$ contains
no singularity of $\Q^0$, we have that \begin{equation}
\label{eq:new1} \lim_{L_k\to 0} f_B(\Q^{L_k}(\x)) = 0 \quad \x\in
K \end{equation} and the limit is uniform on $K$.

Consider a sequence of uniaxial global Landau-de Gennes energy
minimizers $\Q^{L_k}$ such that $L_k\to 0$ as $k\to \infty$. Then
$\Q^{L_k}$ can be written in the form \begin{equation}
\label{eq:new3}\Q^{L_k} = s^{L_k}\left(\n^{L_k} \otimes \n^{L_k} -
\frac{1}{3}\I \right)\end{equation} for $s^k:\Omega \to \Rr$ and
$n^k\in W^{1,2}(\Omega;S^2)$. Then (\ref{eq:new1}) implies that
(up to subsequence), $s^k$ converges uniformly to $s_+$ everywhere
away from $S_0$ i.e. we have
\begin{equation}
\label{eq:scalar} \left| s^k(\x) - s_+ \right|\leq \eps(L_k,\x)
\quad \x \in \bar{\Omega}\setminus B_\delta (\bf{S_0})
\end{equation}
where $ \eps \to 0^+$ as $k\to \infty$, $B_\delta(\bf{S_0})$ is a
small $\delta$-neighbourhood of the singular set $S_0$ and
$0<\delta<1$ is an arbitrary small constant independent of $L$.

\noindent \textbf{Uniform convergence in the interior
\cite{amaz}:} Let $\left\{\Q^{L_k}\right\}$ be a sequence of
global Landau-de Gennes minimizers in $\Acal_\Q$ such that $L_k\to
0$ as $k\to \infty$. Then (up to a subsequence) $\Q^{L_k} \to
\Q^0$ strongly in $W^{1,2}(\Omega,\bar{S})$.

Let $K\subset\Omega$ be a compact set which does not contain any
singularities of $\Q^0$. We define $$e_L(\Q^L) =
\frac{1}{2}|\grad\Q|^2 + \frac{f_B(\Q)}{L}.$$ Then
 \newline \noindent (i) there exists a constant $C>0$ independent
 of $L$ such that
 \begin{equation}
 \label{eq:bochner}
 -\Delta e_L(\Q^L)(\x) \leq C e_L^2\left(\Q^L\right)(\x) \quad
 \x\in K
 \end{equation}
 for $L$ sufficiently small;
 \newline \noindent (ii) we have a uniform bound for $e_L(\Q^L)$
 in the interior of $\Omega$, away from $S_0$ in the limit $L\to 0^+$ i.e.
 \begin{equation}
\label{eq:uniform3} e_L(\Q^L)(\x)  \leq C^{'}(a^2,b^2,c^2,\Omega)
\quad \x\in K
\end{equation}
for all $L$ sufficiently small and a positive constant $C^{'}$
independent of $L$;
\newline \noindent(iii) $\Q^{L_k}$ converges uniformly to $\Q^0$
everywhere in the interior of $\Omega$, away from $S_0$.
\begin{equation}
\label{eq:uniform2} \lim_{k\to \infty} \Q^{L_k}(\x) = \Q^0(\x) ~
uniformly~ for~ \x\in K.
\end{equation} We emphasize that (\ref{eq:uniform3}) and (\ref{eq:uniform2}) only
hold in the interior of $\Omega$. In Section~\ref{sec:results2},
we extend these uniform convergence results up to the boundary for
the uniaxial case.

A uniaxial global Landau-de Gennes minimizer $\Q^L$ is fully
characterized by its scalar order parameter $s^L$ and
distinguished eigenvector $\n^L$. The scalar order parameter,
$s^L$, is a locally Lipschitz function of $\Q^L$ and hence, is
continuous on $\bar{\Omega}$ \cite{sun}. From \cite{bm}, we have
that
\begin{equation}
\label{eq:new4} 0\leq s^L(\x) \leq s_+ \quad \x\in\bar{\Omega}
\end{equation}
and let $\Theta_L = \left\{\x\in\Omega;s^L(\x)=0\right\}$ denote
the isotropic set of $\Q^L$. We have a topologically non-trivial
boundary condition $\Q_b$ in (\ref{eq:9}) and hence, every
interior extension of $\Q_b$ must have discontinuities. We
interpret the defect set of $\Q^L$ as being the defect set of
$\n^L$. Let $S_{\n}^L$ denote the defect set of $\n^L$ and let
$\x_\n\in S_{\n}^L$. Then $\Q^L(\x_\n) =0$, since $\Q^L$ is
well-defined on $\bar{\Omega}$ and consequently $s^L(\x_\n) = 0$.
We deduce that $S_{\n}^L\subset \Theta_L$ and from \cite{nomizu},
we have that $\n^L$ is analytic everywhere away from $\Theta_L$.
We first make an elementary observation about the defect locations
as $L\to 0^+$.

\begin{lemma}
\label{lem:uniaxialdefect} Let $S_0$ denote the singular set of
the limiting harmonic map $\Q^0$ defined in (\ref{eq:Q0}). Let
$\x_\n \in S_{\n}^L$. Then $$dist\left(\x_\n, S_0\right)\leq
\eps(L)$$ where $\eps(L)\to 0$ as $L\to 0^+$.
\end{lemma}

\textit{Proof:} Let $\x_\n \in S_{\n}^L$. As mentioned above,
$s^L(\x_\n)=0$ and $\x_\n \in\Theta_L$, where $\Theta_L$ has been
defined above. However, the bulk energy density $f_B(\Q^L)$
converges uniformly to its minimum value, everywhere away from
$S_0$, in the interior and up to the boundary, as $L\to 0$.
Recalling (\ref{eq:scalar}), we deduce that $dist(\x_\n, S_0) \to
0$ as $L\to 0^+$. Lemma~\ref{lem:uniaxialdefect} now follows.
$\Box$

Lemma~\ref{lem:uniaxialdefect} is also equivalent to the statement
$dist(\Theta_L, S_0) \to 0$ as $L\to 0^+$ i.e. the isotropic set
of a uniaxial global Landau-de Gennes energy minimizer converges
to the singular set of a limiting harmonic map in the limit of
vanishing elastic constant.

\begin{proposition}
\label{prop:1}Let $\Q^L$ be an uniaxial global minimizer of
$\Ical$ in the admissible space $\Acal_\Q$, for a fixed $L>0$.
Then $\Q^L = s^L\left(\n^L\otimes \n^L - \frac{1}{3}\I\right)$ for
some non-negative function $s:\overline{\Omega}\to \Rr^+$ and
$\n\in W^{1,2}(\Omega;S^2)$. The following equations hold
everywhere in $\Omega$, away from the isotropic set $\Theta_L$:
\begin{eqnarray}
\label{eq:uneq} && \Delta s^L - 3s^L |\grad \n^L|^2 =
\frac{s^L}{3L}\left( 2c^2 (s^L)^2 - b^2 s^L - 3a^2\right)
\\
&& \Delta \n^L_j + |\grad \n^L|^2\n^L_j + 2\frac{\partial_k
s^L}{s^L}\n^L_{j,k} = 0 \quad j,k=1,2,3 \label{eq:uneq1}.
\end{eqnarray} Here $\n^L_{j,k}$ denotes the partial derivative
$\frac{\partial \n^L_j}{\partial \x_k}$. Alternatively, $\n^L =
\left(\sin \theta^L \cos\phi^L, \sin\theta^L
\sin\phi^L,\cos\theta^L\right)$, where $\theta^L,\phi^L$ are
functions of spherical polar coordinates $(r,\theta,\phi)$
centered at the origin. Then $\theta^L$ and $\phi^L$ satisfy the
following coupled nonlinear partial differential equations
\begin{eqnarray}
\label{eq:uneq3} && \grad \cdot \left((s^L)^2 \grad
\theta^L\right) = (s^L)^2 \sin \theta^L \cos \theta^L|\grad \phi^L
|^2 \\ && \grad \cdot \left((s^L)^2 \sin^2 \theta^L \grad \phi^L
\right) =
0. 
\label{eq:uneq4}
\end{eqnarray}
\end{proposition}

\textit{Remark: In general, $\Q \in W^{1,2}$ implies that the
tensor $\n\otimes \n \in W^{1,2}$. However, for simply-connected
three-dimensional domains, $\n\otimes \n \in
W^{1,2}\left(\Omega\right) \implies \n \in
W^{1,2}\left(\Omega;S^2\right)$ \cite{bz}.}

\textit{Proof:} In what follows, we drop the superscript $L$ from
$\Q^L$ for brevity. Since $\Q$ is a classical solution of
(\ref{eq:decoupled}), we have that
\begin{eqnarray}
\label{eq:deriv} && \Q_{ij,k} = \partial_k s\left(\n_i\n_j -
\frac{1}{3}\delta_{ij}\right) + s\left(\n_i\n_{j,k} +
\n_j\n_{i,k}\right) \nonumber
\\ && \Q_{ij,kk} = \Delta s \left(\n_i\n_j -
\frac{1}{3}\delta_{ij}\right) + 2\partial_k s \left(\n_i\n_{j,k} +
\n_j\n_{i,k}\right) + s\left(\n_i\n_{j,kk} + \n_j\n_{i,kk} +
2\n_{i,k}\n_{j,k}\right)
\end{eqnarray} where $ i,j,k=1\ldots 3, \Q_{ij,k} = \frac{\partial \Q_{ij}}{\partial
\x_k}$ etc.

Consider the decoupled equations (\ref{eq:decoupled})
$$L\Q_{ij,kk} = \frac{1}{3}\left(2c^2 s^2 -
b^2s - 3a^2\right)\Q_{ij}$$ and multiply both sides by $\n_i$ to
get the following vector equation
\begin{equation}
\label{eq:deriv10} \frac{2}{3}\n_j\Delta s + 2\partial_k s
\n_{j,k} + s\left(\n_{j,kk} - |\grad\n|^2 \n_j \right) =
\frac{2s}{9L}\left(2c^2 s^2 - b^2 s - 3a^2\right)\n_j.
\end{equation}
Multiplying both sides of (\ref{eq:deriv10}) by $\n_j$, we obtain
the following scalar equation for $s$:-
\begin{equation}
\label{eq:deriv3} \frac{2}{3}\Delta s - 2s|\grad \n|^2 = \frac{2
s}{9L}\left(2c^2 s^2 - b^2s - 3a^2\right)
\end{equation} and (\ref{eq:uneq}) now follows. In
(\ref{eq:deriv10}) and (\ref{eq:deriv3}), we use (\ref{eq:deriv})
and the relations $\n_i\n_i=1, \n_i\n_{i,k}=0$ and $\n_i
\n_{i,kk}=-|\grad \n|^2$.

For (\ref{eq:uneq1}), we multiply both sides of the vector
equation (\ref{eq:deriv10}) by the derivative $\n_{j,p}$ for
$p=1,2,3$ to get the following system of three equations -
\begin{equation}
\label{eq:deriv5} 2 \partial_k s \n_{j,p}\n_{j,k} + s
\n_{j,p}\n_{j,kk} =0 \quad p=1,2,3.
\end{equation}
Multiplying both sides by the scalar order
parameter $s$, (\ref{eq:deriv5}) simplifies to
\begin{equation}
\label{eq:deriv6} \n_{j,p} \partial_k\left(s^2 \n_{j,k}\right) = 0
\quad p=1,2,3.
\end{equation}

Next, we note that for a fixed $p$, $(\n_j, \n_{j,p}, \e_j)$ form
an orthogonal basis at each point $\x\in \Omega$, away from the
isotropic set $\Theta_L$ so that
\begin{equation}
\label{eq:new10}
\partial_k\left(s^2 \n_{j,k}\right) = \la_1 \n_j + \la_2 \e_j
\end{equation}
where $$\la_1 = \n_j\partial_k\left(s^2 \n_{j,k}\right) = -s^2
|\grad \n|^2.$$ We substitute (\ref{eq:new10}) into
(\ref{eq:deriv10}) to get
$$\frac{2s}{3}\n_j\Delta s -2 s^2
|\grad \n|^2 \n_j + \la_2 \e_j   = \frac{2s}{9L}\left(2c^2 s^2 -
b^2 s - 3a^2\right)\n_j$$ from which we deduce that $\la_2 = 0$.
Hence
\begin{equation}
\label{eq:new11}
\partial_k\left(s^2 \n_{j,k}\right) + s^2
|\grad \n|^2 \n_j = 0 \quad j=1\ldots 3
\end{equation}
from which (\ref{eq:uneq1}) follows.

An alternative formulation of (\ref{eq:deriv6}) can be obtained by
writing the unit-vector field $\n$ in terms of its spherical
angles, $\theta^L(r,\theta,\phi)$ and $\phi^L(r,\theta,\phi)$,
where $(r,\theta,\phi)$ are spherical polar coordinates centered
at the origin i.e.
\begin{equation}
\label{eq:deriv8} \n = \left(\sin\theta^L\cos\phi^L,
\sin\theta^L\sin\phi^L, \cos\theta^L\right).
\end{equation}
Straightforward computations show that
\begin{eqnarray}
\label{eq:deriv9} && \frac{\partial \n}{\partial \x_k} =
\partial_k \theta^L\left(\cos\theta^L\cos\phi^L,
\cos\theta^L\sin\phi^L, -\sin\theta^L\right) +
\sin\theta^L~\partial_k \phi^L\left(-\sin\phi^L, \cos\phi^L,
0\right) \nonumber \\ && \frac{\partial^2 \n}{\partial
\x_k\partial \x_k} =
\partial_{kk}\theta^L\left(\cos\theta^L\cos\phi^L,
\cos\theta^L\sin\phi^L, -\sin\theta^L\right) -
\left(\partial_k\theta^L\right)^2\left(\sin\theta^L\cos\phi^L,
\sin\theta^L\sin\phi^L, \cos\theta^L\right)+ \nonumber \\
&& +
2\cos\theta^L\partial_k\theta^L\partial_k\phi^L\left(-\sin\phi^L,
\cos\phi^L, 0\right)+
\sin\theta^L\partial_{kk}\phi^L\left(-\sin\phi^L, \cos\phi^L,
0\right) -
\sin\theta^L\left(\partial_k\phi^L\right)^2\left(\cos\phi^L,
\sin\phi^L, 0\right).
\end{eqnarray}

Substituting (\ref{eq:deriv9}) into (\ref{eq:deriv10}) and taking
the dot product of both sides with $\left(\cos\theta^L\cos\phi^L,
\cos\theta^L\sin\phi^L, -\sin\theta^L\right)$, we obtain
\begin{equation}
\label{eq:deriv11} 2\partial_k s \partial_k \theta^L + s
\partial_{kk}\theta^L - s
\sin\theta^L\cos\theta^L|\grad\phi^L|^2 = 0.
\end{equation} We multiply both sides of (\ref{eq:deriv11}) by
$s$ and equation (\ref{eq:uneq3}) now follows. Similarly, we take
the scalar product of (\ref{eq:deriv10}) with the unit-vector
$\left(-\sin\phi^L, \cos\phi^L, 0\right)$ to obtain
\begin{equation}
\label{eq:deriv12} s \sin\theta^L
\partial_{kk}\phi^L + 2 s \cos\theta^L
\partial_k \theta^L\partial_k\phi^L+2
\sin\theta^L\partial_k s \partial_k\phi^L = 0.
\end{equation}
As above, we multiply both sides of (\ref{eq:deriv12}) by $s \sin
\theta$ and (\ref{eq:uneq4}) then follows. The proof of
Proposition~\ref{prop:1} is now complete. $\Box$

\textit{Comment: We note that for $s$ constant, (\ref{eq:uneq1})
is equivalent to the harmonic map equations $\Delta
\n_0 + |\grad \n_0|^2\n_0 = 0$ \cite{bcl}.  }

Next, we use asymptotic methods to predict the minimizer profile
near isolated isotropic points in $\Theta_L$ and establish a $1-1$
correspondence between isolated isotropic points and isolated
point defects.

\begin{proposition}
\label{prop:2}Let $\Q^L$ be an uniaxial global minimizer of
$\Ical$ in the admissible space $\Acal_\Q$, for a fixed $L>0$.
Then $\Q^L(\x) = s^L(\x)\left(\n^L(\x)\otimes \n^L(\x) -
\frac{1}{3}\mathbf{I}\right)$. Let $\Gamma_\n \subset \Omega$
denote the set of isolated point defects in $\n^L$ and let $\Gamma_L
\subset \Omega$ denote the set of isolated isotropic points of
$\Q^L$.

\noindent (i) Let $\x_L \in \Gamma_L$ be an isolated interior
isotropic point.
 Let $(r,\theta,\phi)$ denote a local spherical co-ordinate system
centered at $\x_L$;, then
\begin{equation}
\label{eq:sing1} \left|\grad \n^L \right|^2 \sim
\frac{\alpha(\theta,\phi)}{r^2}~ \quad \textrm{as $r\to 0$}
\end{equation}
where $\alpha$ only depends on $\theta,\phi$ and is independent of
the radial coordinate $r$. Then $\x_L \in \Gamma_\n$ too.
\newline \noindent (ii) Let $\x_\n \in \Gamma_\n$ be an isolated
point defect. Then $\x_\n \in \Gamma_L$ and hence, $\Gamma_\n =
\Gamma_L$.
\end{proposition}

\textit{Proof:} (i) Consider the coupled equation (\ref{eq:uneq}):
\begin{equation}
\label{eq:sing2} \Delta s - 3s |\grad \n|^2 = \frac{s}{3L}\left(
2c^2 s^2 - b^2 s - 3a^2\right)
\end{equation}
Let $\x_L \in \Gamma_L$ be an isolated isotropic point. Since
$s^L(\x) = \frac{3}{2}\Q^L(\x)\left(\n^L(\x)\otimes \n^L(\x) -
\frac{1}{3}\mathbf{I}\right)$ is the product of two analytic
matrices away from $\x_L$, we deduce that $s^L(\x)$ is analytic
for $0< r<r_0$, for some $r_0>0$. We are interested in the
leading-order behaviour of $|\grad \n^L|^2$ as $r\to 0$.

From the local analyticity of $s^L$, we have the following power
series expansion
\begin{equation}
\label{eq:sing3} s^L(\x) = r^n g(\theta,\phi) +
h\left(r,\theta,\phi\right) \quad 0<r<r_b<r_0,~ n\geq 1,
\end{equation} where $(r,\theta,\phi)$ is a local spherical
coordinate system centered at $\x_L$ , $\left|\frac{h}{r^n
g}\right| = o(1)$ as $r\to 0$ and $r_b$ is the radius of
convergence of the series (\ref{eq:sing3}). Further, $g\neq 0$ for
$r\neq 0$ and $\left\{g,h\right\}$ are analytic functions with
bounded derivatives in a sufficiently small neighbourhood of
$\x_L$. Substituting (\ref{eq:sing3}) into (\ref{eq:sing2}) and
expressing $\Delta s^L$ in spherical polar coordinates
\begin{equation}
\label{eq:deltas} \Delta s^L = \frac{1}{r^2}\frac{\partial}{\partial r}\left(r^2
\frac{\partial s^L}{\partial r}\right) +
\frac{1}{r^2\sin^2\phi}\frac{\partial^2 s^L}{\partial \theta^2} +
\frac{1}{r^2\sin\phi}\frac{\partial}{\partial \phi}\left(\sin\phi
\frac{\partial s^L}{\partial\phi}\right),
\end{equation} we have that
\begin{eqnarray}
\label{eq:sing4}
&& r^{n-2}\left[ n(n+1)g + \frac{g_{_{\theta\theta}}}{\sin^2\phi} + g_{\phi\phi} + \cot \phi g_{\phi} \right] + \nonumber \\
&& + \frac{2 h_r}{r} + h_{rr} + \frac{h_{\theta\theta}}{r^2 \sin^2 \phi} + \frac{h_{\phi\phi}}{r^2} + \cot\phi \frac{h_\phi}{r^2} - \nonumber \\
&& - 3|\grad \n^L|^2 \left(r^n g(\theta,\phi) + h(r,\theta,\phi) \right) = \nonumber\\
&& = \frac{\left(r^n g + h\right)}{3L}\left[ 2c^2\left(r^n~g +
h\right)^2 - b^2\left(r^n~g + h\right) - 3a^2\right] ~as~r\to 0.
\end{eqnarray}

All the terms on the right-hand side are $O(r^n)$ whereas the
leading order term on the left-hand side of (\ref{eq:sing4}) is
$O(r^{n-2})$. Since $h$ is an analytic function and
$\left|\frac{h}{r^n}\right| = o(1)$ as $r\to 0$, we have that
$h_r/r, h_{rr} = o(r^{n-2})$ as $r\to 0$. Therefore, for
(\ref{eq:sing4}) to hold as $r\to 0$, we must have
\begin{equation}
\label{eq:sing5} |\grad \n^L|^2 \sim \frac{1}{3r^2}\left[n(n+1) +
\frac{g_{_{\theta\theta}}}{g \sin^2\phi} + \frac{g_{\phi\phi}}{g}
+ \cot \phi \frac{g_{\phi}}{g} \right]~as~r\to 0
\end{equation} and
(\ref{eq:sing1}) now follows. It follows that $|\grad \n^L|^2$ is
not defined as $r\to 0$ and hence, the isolated isotropic point
$\x_L \in \Gamma_\n$ too.
\newline (ii) Let $\x_\n \in \Gamma_\n$. Then $\Q^L(\x_\n) = 0$,
since $\Q^L$ is well-defined on $\bar{\Omega}$. Therefore, we must
have $s^L(\x_\n) = 0$ and by definition, $\x_\n \in \Gamma_L$.
Combining (i) and (ii), we conclude that $\Gamma_L = \Gamma_\n$.
$\Box$

\noindent\textit{Comment: By analogy with \cite{millot}, one might
expect that uniaxial global minimizers can only account for
isolated point defects and all higher-dimensional defects are
intrinsically biaxial. We hypothesize that $n=2$ in
(\ref{eq:sing3}) i.e. we have a quadratic decay of the scalar
order parameter as we approach point defects, by analogy with the
study of vortices in Ginzburg-Landau theory \cite{maj2}.}

\noindent \textit{Comment: The estimate (\ref{eq:sing1}) is
analogous to a similar result on singularity profiles within the
Oseen-Frank theory of uniaxial nematic liquid crystals with a
constant scalar order parameter $s$ \cite{bcl}. In \cite{bcl}, the authors
show that near every singularity $\x_p \in \Omega$, we have
\begin{equation}
\label{eq:hedgehog}\n \sim R\frac{\x - \x_p}{|\x - \x_p|}\end{equation}
for some rotation $R\in SO(3)$. Therefore ,
$$|\grad \n|^2 \sim \frac{2}{|\x - \x_p|^2}~as~\x \to \x_p.$$
The estimate (\ref{eq:sing1}) suggests that we have a similar
radial hedgehog-type of profile (\ref{eq:hedgehog})  near the isolated zero $\x_L \in
\Omega$, for uniaxial global minimizers within the Landau-de
Gennes theory.}

\section{Far-field results}
\label{sec:results2}

In this section, we study the qualitative properties of uniaxial
global minimizers $\left\{\Q^L\right\}$ away from the isotropic set
$\Theta_L$, in the limit $L\to 0$. This is equivalent to studying
the qualitative properties of $\left\{\Q^L\right\}$ away from the
singular set $S_0$ of the limiting harmonic map $\Q^0$ defined in
(\ref{eq:Q0}), as $L\to 0^+$.

 Let $\Q^L = s^L \left(\n^L\otimes \n^L - \frac{1}{3}\mathbf{I}\right)$ be a uniaxial global minimizer
 for fixed $L>0$. Recall that for $L$
sufficiently small,
\begin{equation}
\label{eq:rsnew} 0 \leq s_+ - s^L(\x) \leq \eps_1(L)
\end{equation}
or equivalently
\begin{equation}
\left||\Q^L|^2 - \frac{2}{3}s_+^2 \right| \leq \eps_2(L)
\label{eq:rs1}
\end{equation}
where $\eps_1(L),\eps_2(L) \to 0$ as $L\to 0^+$, everywhere away
from $S_0$.

Our first result is an inequality for
$$A^L = \frac{1}{2}\Q^L_{ij,k}\Q^L_{ij,k}$$
that holds everywhere away from $S_0$ on $\bar{\Omega}$. We do not
use Lemma~\ref{lem:3} in the subsequent sections but keep it as an
interesting technical result.

\begin{lemma}
\label{lem:3} Let $A^L = \frac{1}{2}\Q^L_{ij,k}\Q^L_{ij,k}$ by
definition. Then we have the
following inequality on $\Omega\setminus B_\delta(S_0)$ for $L$ sufficiently small
\begin{equation}
\label{eq:rs2} -\Delta A^L + \left|D^2 \Q^L\right|^2 \leq
\frac{1}{\alpha^2}\left|D^2 \Q^L\right|^2  + \alpha^4
\frac{A{^L}^2}{|\Q^L|^2}
\end{equation}
where $\alpha>1$ is a positive constant independent of $L$ that
can be worked out explicitly, $B_{\delta}(S_0)$ is a small
$\delta$-neighbourhood of $S_0$ and $\delta>0$ is independent of
$L$.
\end{lemma}

\textit{Proof:} The derivation of (\ref{eq:rs2}) closely follows
the methods in \cite{bbh}. In what follows, we drop the
superscript $L$ for brevity. First, consider the decoupled
equations (\ref{eq:decoupled}); setting
$$f(s) = \left(2c^2s^2 - b^2 s - 3a^2\right)$$ and differentiating
both sides of (\ref{eq:decoupled}) with respect to $\x_p$, we
obtain
\begin{equation}
\label{eq:rs3} \Q_{ij,kkp} = \frac{\Q_{ij,p}}{3L}f(s) +f^{'}(s)
\frac{\Q_{ij}}{\sqrt{6}L}
\frac{\Q_{rs}\Q_{rs,p}}{|\Q|}~for~p=1,2,3.
\end{equation}
From (\ref{eq:rs1}) and the global upper bound (\ref{eq:max}), we
have that $|\Q|$ is bounded away from zero on $\Omega \setminus
B_\delta(S_0)$ and
\begin{eqnarray}
\label{eq:rs4}
 && f(s)\leq 0  \nonumber
\\
&& f^{'}(s) > 0  \nonumber
\\ && f^{''}(s) > 0,
\end{eqnarray} on the set $\Omega \setminus B_\delta(S_0)$, where $f^{'}(s) = \frac{df}{ds}$, $f^{''}(s) =
\frac{d^2 f}{ds^2}$ etc. A straightforward computation shows that
\begin{equation}
\label{eq:rs5} \Delta A = |D^2 \Q|^2 + \Q_{ij,k}\Q_{ij,ppk}
\end{equation}
and using (\ref{eq:rs3}), we obtain
\begin{equation}
\label{eq:rs6} \Delta A =|D^2 \Q|^2 +  |\grad \Q|^2\frac{f(s)}{3L}
+
\frac{f^{'}(s)}{\sqrt{6}L}\frac{\left(\Q\cdot\grad\Q\right)^2}{|\Q|}.
\end{equation}
From (\ref{eq:rs4}) and (\ref{eq:decoupled}), we have the
following inequality
\begin{equation}
\label{eq:rs7} -\Delta A + |D^2 \Q|^2 \leq |\grad \Q|^2
\frac{|\Delta \Q|}{|\Q|}.
\end{equation}
Finally, we use the inequality
$$|\Delta \Q| \leq \alpha \left|D^2 \Q \right|$$
where $\alpha>1$ is a positive constant that can be worked out
explicitly. Substituting the above into (\ref{eq:rs7}),
\begin{equation}
\label{eq:rs8} -\Delta A + |D^2 \Q|^2 \leq 2\alpha A \frac{|D^2
\Q|}{|\Q|} \leq \frac{1}{\alpha^2}|D^2 \Q|^2 + \alpha^4
\frac{A^2}{|\Q|^2}
\end{equation}
and (\ref{eq:rs2}) now follows. $\Box$

We recall the uniform convergence result in (\ref{eq:uniform2})
and (\ref{eq:uniform3}), whereby we establish a uniform bound for
$|\grad \Q^L|$, independent of $L$, everywhere away from $S_0$ in
the interior of $\Omega$. The next step is to extend this uniform
convergence result up to the boundary. To do so, we adapt the
small energy regularity theorem in \cite{chen} to the Landau-de
Gennes framework to obtain a uniform bound for $|\grad \Q^L|$
independent of $L$, everywhere away from $S_0$ up to the boundary.

Consider a boundary point $\x_0 \in \partial\Omega$ and define the
region $\Omega_r(\x_0) = \bar{\Omega}\cap B_r(\x_0)$, where
$B_r(\x_0)$ is a ball of radius $r$ centered at $\x_0$. Let $\rho$
be a suitably small positive constant such that for any
$\x_0\in\partial\Omega$, we may choose a coordinate system
$\left\{\x_\alpha\right\}$ so that $\x_0$ is at the origin and
$\Omega_\rho(\x_0)$ corresponds to $B_\rho^+(\x_0) =
\left\{\x\in\bar{\Omega}; |\x|\leq \rho;~ x_3 \geq 0 \right\}$.
\begin{proposition}
\label{prop:chen} Let $\left\{\Q^{L_k}\right\}$ be a sequence of uniaxial global
minimizers for $\Ical$ in the admissible space $\Acal_\Q$, where
$L_k \to 0$ as $k\to \infty$. We can extract a subsequence such
that $\Q^{L_k} \to \Q^0$ strongly in $W^{1,2}(\Omega;\bar{S})$ as
$k\to\infty$. Let $\x_0\in\partial\Omega$ be such that
$\Omega_r(\x_0)$ contains no singularity of the limiting harmonic
map $\Q^0$. Then there exists $C_1>0, C_2>0, r_0>0, \bar{L}_0
>0$ (all constants independent of $L_k$) so that if
\begin{equation}
\label{eq:chen1} \int_{\Omega_r(\x_0)}\frac{1}{2}|\grad \Q^L|^2 +
\frac{f_B(\Q^L)}{L}~d\x \leq C_1 \quad r<\min\left\{r_0,\rho\right\}
\end{equation} then
\begin{equation}
\label{eq:chen2} r^2 \sup_{\Omega_{r/2}(\x_0)}e_L(\Q^L) \leq C_2
\quad~for~all~L_k< \bar{L}_0
\end{equation}
where $$e_L(\Q^L) = \frac{1}{2}|\grad \Q^L|^2 +
\frac{f_B(\Q^L)}{L}.$$
\end{proposition}

\noindent{\textit{Proof:}} The first half of the proof of
Proposition~\ref{prop:chen} closely follows the scaling arguments
for the interior uniform convergence result (\ref{eq:uniform3}) in
\cite{amaz} and the second half closely follows the arguments in
Theorem~$2.1$ in \cite{chen}.

We first recall from (\ref{eq:new1}) that since $\Omega_r(\x_0)$
contains no singularity of $\Q^0$, $\exists \m(\x)\in S^2$ such
that
\begin{equation}
\label{eq:chen3} \left|\Q^L(\x) - s_+\left(\m(\x)\otimes \m(\x) -
\frac{1}{3}\mathbf{I}\right)\right| < \eps_0 <<1 \quad
\x\in\Omega_r(\x_0)
\end{equation}
for $L$ sufficiently small.

We continue reasoning similarly to \cite{amaz}. We fix an
arbitrary $L_k<\bar{L}_0$. We let $0<r_1<\frac{2r}{3}< \min\left\{\frac{2r_0}{3}, \frac{2\rho}{3}\right\}$ and $\x_1 \in
\Omega_{r_1}(\x_0)$ be such that
\begin{eqnarray}
\label{eq:chen4} && \max_{0\leq s\leq
\frac{2r}{3}}\left(\frac{2r}{3} - s\right)^2 \max_{\Omega_s(\x_0)}
e_{L_k}\left(\Q^{L_k}\right) = \left(\frac{2r}{3} - r_1 \right)^2
e_{L_k}\left(\Q^{L_k}\right)(\x_1).
\end{eqnarray}
Define $e_{1}^{L_k} = \max_{\x\in\Omega_{r_1}(\x_0)}
e_{L_k}(\Q^{L_k}) = e_{L_k}(\Q^{L_k})(\x_1)$. Then
\begin{eqnarray}
\label{eq:chen5} \max_{\Omega_{\frac{2/3r - r_1}{2}}(\x_1)}
e_{L_k}(\Q^{L_k}) \leq 4 e_{1}^{L_k}
\end{eqnarray}
where we use the inclusion $\Omega_{\frac{2/3r -
r_1}{2}}(\x_1)\subset \Omega_{\frac{2/3r + r_1}{2}}(\x_0)$,
$\frac{2/3r + r_1}{2} \leq \frac{2r}{3}$ by definition of $r_1$
and the inequalities (\ref{eq:chen4}).

Define $r_2 =\frac{2/3r - r_1}{2}\sqrt{e_{1}^{L_k}}$ and let
\begin{equation}
\label{eq:chen6} \R^{L_k}(\x) = \Q^{L_k}\left(\x_1 +
\frac{\x}{\sqrt{e_{1}^{L_k}}}\right).
\end{equation}
Let $\bar{L}_k = e_{1}^{L_k} L_k$. Then $\R^{L_k}$ has the
following properties on $\Omega_{r_2}(0)$:-
\begin{eqnarray}
 && e_{\bar{L}_k}(\R^{L_k}) =
\frac{1}{e_{1}^{L_k}}e_{L_k}(\Q^{L_k}) \label{eq:chen7} \\ &&
\max_{\x\in\Omega_{r_2}(0)}e_{\bar{L}_k}(\R^{L_k}) \leq 4 \quad
e_{\bar{L}_k}(\R^{L_k})(0)=1 \label{eq:chen8} \\ &&
\R^{L_k}_{ij,kk} = \frac{1}{\bar{L}_k}\left(2c^2 s^2 - b^2 s -
3a^2\right) \R^{L_k}_{ij} \label{eq:chen9}
\end{eqnarray}
where $s^2 = \frac{3}{2}\left|\Q^{L_k}\right|^2$.

We next claim that $r_2\leq 1$. It is obvious that $r_2\leq 1$
implies the conclusion (\ref{eq:chen2}). We prove this claim by
contradiction. Assume that $r_2>1$; then using the same arguments
as in \cite{chen}, one is led to the existence of a sequence of
solutions $\left\{\R^{L_k}\right\}$ of (\ref{eq:chen9}) on
$\Omega_{1}(0) = B_1^+(0)$, with the following properties: -
\begin{eqnarray}\label{eq:chen10}
&& - \Delta \R^{L_k}_{ij} + \frac{1}{\bar{L}_k}\left(2c^2 s^2 -
b^2 s - 3a^2\right) \R^{L_k}_{ij} = 0 ~in~ B_{1}^+(0) \nonumber \\
&& \max_{\x\in B_{1}^+(0)}e_{\bar{L}_k}(\R^{L_k}) \leq 4 \quad
e_{\bar{L}_k}(\R^{L_k})(0)=1  \nonumber \\ && \R^{L_k}|_{\x_3 = 0}
= \Q_b\left(\x + \frac{\x_1}{\sqrt{e_1^{L_k}}}\right) ~with~ \nonumber \\ && 
|\grad \R^{L_k}||_{\x_3 = 0}\leq \eps_k |\grad
\Q_b|_{L^{\infty}(\partial\Omega)},~ |\grad^2 \R^{L_k}||_{\x_3 =
0}\leq \eps_k^2 |\grad^2 \Q_b|_{L^{\infty}(\partial\Omega)} ~with~\eps_k\to 0~as~ k\to\infty
\nonumber \\ && \int_{B_{1}^+(0)}e_{\bar{L}_k}(\R^{L_k})~d\x \leq
\delta_k \to 0^+ ~as~ k\to \infty.
\end{eqnarray}
From (\ref{eq:bochner}) and (\ref{eq:chen8}), we deduce that
$\R^{L_k}$ satisfies the following Bochner-type inequality on
$B_{1}^+(0)$:-
\begin{equation}
\label{eq:chen11} -\Delta e_{\bar{L}_k}(\R^{L_k})\leq C^{'}
e_{\bar{L}_k}(\R^{L_k}) \quad \x\in B_{1}^+(0)
\end{equation}
where $C^{'}$ is a constant independent of $L_k$.

Next, we write $\R^{L_k}$ explicitly in terms of its scalar order
parameter and leading eigenvector -
\begin{equation}
\label{eq:chen12} \R^{L_k}_{ij} = s_k \left(\n^k_i\n^k_j -
\frac{1}{3}\delta_{ij}\right) \quad \n^k\in
W^{1,2}\left(\Omega;S^2\right)
\end{equation}
where $|\R^{L_k}|^2 = \frac{2}{3}s_k^2$ and
$$ \left| s_k - s_+\right| \leq \frac{s_+}{100}$$
from (\ref{eq:chen3}), for $k$ sufficiently large. From
Proposition~\ref{prop:1}, we have that $s^k$ and $\n^k$ satisfy
the following equations in $B_1^+(0)$:-
\begin{eqnarray}
&& \Delta s_k - 3s_k|\grad \n^k|^2 = \frac{s_k}{3\bar{L}_k}(2c^2
s_k^2 - b^2 s_k - 3a^2) \label{eq:chen13} \\
&& \Delta \n^k_j + |\grad \n^k|^2 \n^k_j + 2
\frac{\partial_{p}s_k}{s_k}\n^k_{j,p} = 0 \label{eq:chen14} \\
&& |\grad \R^{L_k}|^2 = \frac{2}{3}|\grad s_k|^2 +
2s_k^2|\grad\n^k|^2 \leq 4 \label{eq:chen15}.
\end{eqnarray}

From (\ref{eq:chen3}) and (\ref{eq:chen15}), we deduce that
\begin{eqnarray}
\label{eq:chen16} |\grad \n^k| \leq \frac{2}{s_+},~ \left| \frac{2
\grad s_k}{s_k}\right|\leq \frac{6}{s_+}\quad on~ B_{1}^+(0).
\end{eqnarray}
We combine (\ref{eq:chen16}) and  (\ref{eq:chen14}) to deduce that
\begin{equation}
\label{eq:chen17} \sup_{B_{2/3}^+(0)}|\grad \n^k|^2 \leq c
\delta_k^{1/3} \to 0 \quad k\to \infty
\end{equation} where $c$ is a constant independent of $k$. In particular, this implies that $||\grad
\n^k||_{L^{\infty}(\partial\Omega)}\leq c\delta_k^{1/3} \to 0$ as
$k\to\infty$. The proof of (\ref{eq:chen17}) is identical to the
proof of Theorem~$2.1$ in \cite{chen} and the details are omitted
for brevity.

Next, look at the equation (\ref{eq:uneq}) and introduce the
function
$$ \bar{s}_k = s_+ - s_k.$$
Then $\bar{s}_k$ is a solution of the following problem on
$B_{1}^+(0)$ :-
\begin{eqnarray}
\label{eq:chen18} && -\Delta \bar{s}_k = 3s_k |\grad \n^k|^2 -
\frac{2 c^2 s_k \bar{s}_k (s_k - s_{-})}{3\bar{L}_k} \nonumber\\
&& \bar{s}_k(\x) = 0 \quad \x\in\left\{B_1^{+}(0)\cap x_3 =
0\right\}
\end{eqnarray} where $s_{-}<0$ is a constant.
Repeating the same arguments as in \cite{chen}, we obtain the
following estimates:-
\begin{eqnarray}
\label{eq:chen19} && \bar{s}_k(\x)\leq c_1 x_3 \delta_k^{1/3}
\quad \x\in B_{1/2}^+(0) \nonumber \\ && ||\grad
\bar{s}_k||_{L^{\infty}(\partial\Omega)} = ||\grad
s_k||_{L^{\infty}(\partial\Omega)} \leq c_0 \delta_k^{1/3} \to 0
\quad k\to\infty
\end{eqnarray} where $c_0$ and $c_1$ are positive constants
independent of $L$.

Finally, we define $\tilde{e}_k = \max\left\{0, e_{L_k}(\R^{L_k})
- (c + c_0) \delta_k^{1/3}\right\}$. (Recall that $f_B(\R^{L_k})=0$ on $x_3=0$ because of the choice of the boundary condition $\Q_b$.) Then from
(\ref{eq:bochner}), (\ref{eq:chen17}) and (\ref{eq:chen19}), we have that
\begin{eqnarray}
&& -\Delta \tilde{e}_k(\x) \leq C^{''} \tilde{e}_k(\x) \quad \x\in
B_{1/2}^+(0)\nonumber \\ && \tilde{e}_k(\x)|_{x_3=0} = 0
\label{eq:chen20}
\end{eqnarray} for a constant $C^{''}$ independent of $L$.
Using standard arguments as in \cite{chen}, (\ref{eq:chen20})
implies that
$$\sup_{\x\in B_{1/4}^+(0)}\tilde{e}_k(\x) \leq c_3 \delta_k \to
0^+~as ~ k\to \infty$$ contradicting (\ref{eq:chen8}). The proof
of Proposition~\ref{prop:chen} is now complete. $\Box$

For the reader's convenience, we quote Lemma~$2$ from \cite{bbh}
which is used in Proposition~\ref{prop:convergence}:

\textit{Lemma $2$ from \cite{bbh}: Let $\omega(r)$ be a solution of
\begin{eqnarray}
&& -\eps^2 \Delta\omega + \omega = 0 ~on~ B(0,R) \nonumber \\ && \omega = 1 ~on~ \partial B(0,R)
\label{eq:bbh1}.
\end{eqnarray} Then for $\eps<\frac{3R}{4},
\omega(r) \leq e^{\frac{1}{4\eps R}\left(r^2 - R^2\right)}$ on
$B(0,R)$.}

\begin{proposition}
\label{prop:convergence} Let $\left\{\Q^{L_k}\right\}$ be a
sequence of uniaxial global Landau-de Gennes minimizers in the
admissible space $\Acal_\Q$, where $L_k\to 0^+$ as $k\to \infty$.
Then as $k\to\infty$, we can extract a suitable subsequence such
that $\Q^{L_k} \to \Q^0$ in
$C^{1,\alpha}\left(\overline{\Omega}\setminus
B_\delta(S_0)\right)$ for some $0 <\alpha < 1$ and
$B_{\delta}(S_0)$ is a small $\delta$-neighbourhood of the
singular set, $S_0$, of the limiting harmonic map, $\Q^0$, where
$\Q^0$ has been defined in (\ref{eq:Q0}) and $\delta$ is
independent of $L_k$.
\end{proposition}

\textit{Proof:} The proof follows the methods in \cite{bbh} and the key ingredient is to establish a global bound
for $\frac{s_+ - s^{L_k}}{L_k}$, everywhere away from $S_0$, for
$L_k$ sufficiently small.

We drop the superscript $L_k$ in what follows for convenience.
Consider the equation (\ref{eq:uneq}) on
$\overline{\Omega}\setminus B_{\delta}(S_0)$ and introduce the
function
\begin{equation}
\label{eq:psi} \psi = \frac{s_+ - s}{L}, \end{equation} $s_+$ has
been defined in (\ref{eq:8}) and $s_- = \left(b^2 - \sqrt{b^2 +
24a^2 c^2}\right)/4c^2.$ Then (\ref{eq:uneq}) can be re-written as
\begin{equation}
\label{con1} \Delta s - 3s|\grad \n|^2 = - \frac{2c^2 s}{3} \psi
\left(s - s_{-} \right)
\end{equation}
From (\ref{eq:max}) and (\ref{eq:scalar}), we have that
$\frac{2}{3}s_+^2\geq |\Q|^2 \geq \frac{2}{3}s_+^2 - \eps_L$ where
$\eps_L \to 0$ as $L \to 0$, on $\overline{\Omega}\setminus
B_{\delta}(S_0)$. Therefore,
$$ \frac{2c^2 s}{3}\left(s - s_{-} \right)\geq
\frac{1}{\beta}$$ on $\overline{\Omega}\setminus B_{\delta}(S_0)$,
where $\beta$ is a positive constant independent of $L$.

We note that $|\grad \Q|^2 = \frac{2}{3}|\grad s|^2 + 2s^2|\grad
\n|^2$ where $\Q = s\left(\n\otimes \n -
\frac{1}{3}\mathbf{I}\right)$. We recall the global uniform bound
(\ref{eq:chen2}) everywhere away from $S_0$ to deduce that
$$|\grad \n|^2 \leq C(a^2,b^2,c^2,\Omega) $$
on $\overline{\Omega}\setminus B_\delta(S_0)$. Combining the
above, we have that $\psi$ satisfies the following inequality on
$\overline{\Omega}\setminus B_\delta(S_0)$
\begin{equation}
\label{psi1} -\beta L \Delta \psi + \psi \leq \gamma |\grad \n|^2
\leq D(a^2,b^2,c^2,\Omega)
\end{equation}
where $\gamma$ and $D$ are positive constants independent of $L$.
Applying standard maximum principle arguments, we conclude that
\begin{equation}
\label{eq:psi10} \|\psi\|_{L^{\infty}(\Omega\setminus
B_\delta(S_0))} \leq D^{'}(a^2,b^2,c^2,\Omega)
\end{equation}
where $D^{'}$ is a positive constant independent of $L$.

Consider the governing equations (\ref{eq:decoupled}) for a
uniaxial global minimizer $\Q$; they can be written in terms of
the function $\psi$ as shown below -
\begin{equation}
\Delta \Q = \frac{1}{3L}\left(2c^2s^2 - b^2 s - 3a^2\right)\Q \leq
 - \alpha \psi \Q \label{eq:psi7}
\end{equation}
where $\alpha>0$ is a constant independent of $L$, we have used
the definition of $\psi$ in (\ref{eq:psi}) and the uniform
convergence of bulk energy density everywhere away from $S_0$
(refer to (\ref{eq:new1})). Finally, we combine the global upper
bound (\ref{eq:max}) and the $L^{\infty}$-estimate
(\ref{eq:psi10}) to conclude that
\begin{equation}
\label{eq:psi11} \|\Delta \Q\|_{L^{\infty}(\overline{\Omega}\setminus
B_\delta(S_0))} \leq D^{''}(a^2,b^2,c^2,\Omega)
\end{equation}
where $D^{''}$ is a positive constant independent of $L$ i.e.
$|\Delta \Q|$ can be bounded independently of $L_k$ everywhere
away from $S_0$. Finally, we use (\ref{eq:psi11}) and Sobolev
estimates to establish $\left\{\Q^{L_k}\right\} \to \Q^0$ in
$C^{1,\alpha}(\Omega\setminus B_\delta(S_0))$ as $k\to\infty$, for
some $0<\alpha<1$. The proof of Proposition~\ref{prop:convergence}
is now complete. $\Box$

\textit{Comment: One immediate consequence of (\ref{eq:psi10}) is
that $s_+ - s^L\leq C L$, where $C$ is a positive constant
independent of $L$, everywhere away from $S_0$ in the limit $L\to
0$. This explicitly estimates the rate of convergence in
(\ref{eq:scalar}) and improves upon a previous estimate in
\cite{amaz} where an analysis of the bulk energy density $f_B$ in
(\ref{eq:6}) shows that $s_+ - s\leq C_1\sqrt{L}$ where $C_1$
is a positive constant independent of $L$.}

\begin{lemma}
\label{lem:lem7} Let $\Q^L=s^L\left(\n^L \otimes \n^L -
\frac{1}{3}\I\right)$ be a uniaxial global minimizer of $\Ical$ in
$\Acal_\Q$, for $L$ sufficiently small. Then for $\x \in
\Omega\setminus B_\delta(S_0)$, we have
\begin{eqnarray}
\label{eq:lem7} && |\grad s^L| \leq \eps_1(\x) \nonumber \\
&& \left||\grad \n^L(\x)|^2 - |\grad \n_0|^2 \right| \leq \eps_2(\x)
\end{eqnarray}
where $\n_0$ and $\Q^0$ are defined in (\ref{eq:Q0}) and
$\eps_1,\eps_2 \to 0^+$ as $L\to 0^+$.
\end{lemma}

\textit{Proof:} Lemma~\ref{lem:lem7} is a direct consequence of
Proposition~\ref{prop:convergence}. Let $\x \in \Omega\setminus
B_\delta(S_0)$. Then from Proposition~\ref{prop:convergence}, we
have that
\begin{eqnarray}
\label{eq:lem8} && \left|\Q^L_{ij}(\x)-\Q^0_{ij}(\x)\right|\leq
\eps_3(\x)\nonumber \\ &&
\left|\Q^L_{ij,k}(\x)-\Q^0_{ij,k}(\x)\right|\leq \eps_4(\x)
\end{eqnarray}
where $\Q^0$ is the limiting harmonic map in (\ref{eq:Q0}),
$\Q_{ij,k}=\frac{\partial \Q_{ij}}{\x_k}$ and $\eps_3,\eps_4<<1$.
One can directly compute \begin{equation} \label{eq:7f}\left|\grad
\Q^0\right|^2 = 2s_+^2\left|\grad  \n_0 \right|^2.
\end{equation} On the other hand, $$\left|
\Q^L \right|^2 = \frac{2}{3}(s^L)^2$$ and therefore,
\begin{equation}
\label{eq:lem7a} \Q^L_{ij}\Q^L_{ij,k} = \frac{2}{3} s^L \partial_k
s^L
\end{equation}
where $|s^L(\x) - s_+| < \eps_5(\x) <<1$ for $\x\in
\Omega\setminus B_\delta(S_0)$, from (\ref{eq:scalar}).

One can
re-write $\Q^L_{ij}\Q^L_{ij,k}$ as shown below -
\begin{equation}
\label{eq:lem7b} \Q^L_{ij}\Q^L_{ij,k} =
(\Q^L_{ij}(\x)-\Q^0_{ij}(\x))\Q^L_{ij,k}(\x) +
\Q^0_{ij}(\x)\left(\Q^L_{ij,k}(\x)-\Q^0_{ij,k}(\x)\right)
\end{equation} since $\Q^0_{ij}\Q^0_{ij,k} = 0$ from $|\Q^0|^2 =
\frac{2}{3}s_+^2$. Using the inequalities (\ref{eq:lem8}), the
global bound (\ref{eq:chen2}) and the triangle inequality, we have
that
\begin{equation}
\label{eq:lem7c} \left|\Q^L_{ij}(\x)\Q^L_{ij,k}(\x) \right|\leq
\eps_6(\x) << 1
\end{equation}
for $\x\in \Omega\setminus B_\delta(S_0)$ and from
(\ref{eq:lem7a}) and (\ref{eq:scalar}), this necessarily implies
that
\begin{equation}
\label{eq:7d}|\grad s^L| \leq \eps_7\left( L \right)
\end{equation} away from $S_0$, where $\eps_7 \to 0^+$ as
$L\to 0^+$.

On the other hand, from Proposition~\ref{prop:convergence}, $\Q^L
\to \Q^0$ in $C^{1,\alpha}(\Omega;\bar{S})$ as $L\to 0$ (up to a subsequence), everywhere away
from  $S_0$. Therefore, for $\x\in \Omega\setminus B_\delta(S_0)$,
\begin{eqnarray}
\label{eq:7e} \left| |\grad \Q^L|^2 - |\grad \Q^0|^2\right| \leq
\eps_8 (\x)
\end{eqnarray}
where $\eps_8 \to 0^+$ as $L\to 0^+$. A direct computation shows
that
$$|\grad \Q^L|^2 = \frac{2}{3}|\grad s^L|^2 + 2(s^L)^2|\grad
\n^L|^2.$$ Combining (\ref{eq:scalar}), (\ref{eq:7d}), (\ref{eq:7e}) and
(\ref{eq:7f}), we have that $|\grad \n^L|^2 \to |\grad \n_0|^2$ as
$L\to 0^+$. Lemma~\ref{lem:lem7} now follows. $\Box$

\begin{proposition}
\label{prop:4} Let $\Q^L=s^L\left(\n^L \otimes \n^L -
\frac{1}{3}\I\right)$ be a uniaxial global minimizer of $\Ical$ in
$\Acal_\Q$, for $L$ sufficiently small. Then for $\x \in
\Omega\setminus B_\delta(S_0)$, we have that
\begin{eqnarray}
\label{eq:prop4} \left| \frac{s_+ - s^L}{L} - \frac{9\left|\grad \n_0\right|^2}{\sqrt{b^4
+ 24a^2c^2}}\right| \leq \eps_9(\x)
\end{eqnarray}
where $\eps_9 \to 0^+$ as $L \to 0^+$.
\end{proposition}

\textit{Proof:} Consider the function $\psi^L= \frac{s_+ -
s^L}{L}$ in (\ref{eq:psi}) and the equation (\ref{eq:uneq}) on
$\Omega\setminus B_\delta(S_0)$
\begin{equation}
\label{prop4} \Delta s^L - 3s^L |\grad \n^L|^2 = - 2c^2
\frac{s^L}{3} (s^L - s_{-})\psi
\end{equation} where $s_{\pm} = \frac{b^2 \pm \sqrt{b^4 + 24a^2 c^2}}{4c^2}$.
Equation (\ref{prop4}) can be re-arranged to give
\begin{eqnarray}
\label{prop5} && - L \Delta \left( \psi^L - \frac{9 |\grad
\n_0|^2}{\sqrt{b^4 + 24a^2c^2}}\right) + 2c^2 \frac{s^L(s^L -
s_-)}{3} \left(\psi^L - \frac{9 |\grad \n_0|^2}{\sqrt{b^4 +
24a^2c^2}}\right) \nonumber = \\ &&
 = 3s^L|\grad \n^L|^2 + \frac{9 L}{\sqrt{b^4 +
24a^2c^2}}\Delta \left|\grad \n_0\right|^2 - \frac{6c^2 s^L \left(s^L -
s_-\right)}{\sqrt{b^4 + 24a^2c^2}}\left|\grad \n_0\right|^2.
\end{eqnarray}
We note that $\Delta\left|\grad \n_0\right|^2 = O(1)$ away from
$S_0$ and the right-hand side of (\ref{prop5}) can be written as
\begin{eqnarray}
&& 3s^L|\grad \n^L|^2 + \frac{9 L}{\sqrt{b^4 +
24a^2c^2}}\Delta \left|\grad \n_0\right|^2 - \frac{6c^2 s^L \left(s^L -
s_-\right)}{\sqrt{b^4 + 24a^2c^2}}\left|\grad \n_0\right|^2 = \\ &&  = 3s^L|\grad
\n^L|^2 - 3s_+ \left|\grad \n_0\right|^2 + 3s_+ \left|\grad \n_0\right|^2  - \frac{6c^2 s^L
\left(s^L - s_-\right)}{\sqrt{b^4 + 24a^2c^2}}\left|\grad \n_0\right|^2  +
O(L).\label{prop6}
\end{eqnarray}
Finally, we use (\ref{eq:scalar}) and (\ref{eq:lem7}) to deduce
that
$$ 3s^L|\grad
\n^L|^2 - 3s_+ \left|\grad \n_0\right|^2 \leq \eps_{10}$$
where $\eps_{10}\to 0^+$ as $L\to 0^+$ and
$$ \frac{6c^2 s^L
\left(s^L - s_-\right)}{\sqrt{b^4 + 24a^2c^2}}\left|\grad \n_0\right|^2  \to 3s_+ \left|\grad \n_0\right|^2
$$ as $L\to 0^+$, since
$s^L - s_{-} \to \left(s_+ - s_{-}\right) = \sqrt{b^4 +
24a^2c^2}/2c^2$ as $L\to 0^+$. Combining the above, we have that
\begin{equation}
\label{prop7} -L \Delta\left( \psi^L - \frac{9 |\grad
\n_0|^2}{\sqrt{b^4 + 24a^2c^2}}\right) + \beta \left( \psi^L -
\frac{9 |\grad \n_0|^2}{\sqrt{b^4 + 24a^2c^2}}\right)\leq
\eps_{11}
\end{equation}
where $\beta$ is a positive constant independent of $L$ and
$\eps_{11} \to 0^+$ as $L\to 0^+$. Proposition~\ref{prop:4} now follows from the maximum
principle and Lemma~2 of \cite{bbh}. $\Box$

\section{Generalizations}
\label{sec:dis}

This paper focuses on qualitative properties of global minimizers
of the Landau-de Gennes energy functional on $2D$ and $3D$
domains. In the $2D$ case, we show that the Landau-de Gennes
theory is equivalent to Ginzburg-Landau theory for superconductors
and make predictions about the dimension of the defect set, the
defect locations and the asymptotic profile of global minimizers
close to and far away from the defect set.

In $3D$, we focus on uniaxial global minimizers of the Landau-de
Gennes energy functional because this is the first step in a
rigorous study of arbitrary minimizers. The topological defects
are contained inside the isotropic set of the uniaxial global
minimizers. We derive the governing equations for the scalar order
parameter $s^L$ and the leading eigenvector $\n^L$; these
equations reflect the coupling between the two quantities. We show
that the topological defects (or equivalently the isotropic set) are necessarily contained in a small
neighbourhood of the singular set of a limiting harmonic map and
establish the vortex-like or radial hedgehog-like profile of
isolated point defects. We also study the qualitative properties
of uniaxial global minimizers away from the isotropic set. In
particular, we establish the $C^{1,\alpha}$-convergence of
uniaxial global minimizers to a limiting harmonic map, everywhere
away from the isotropic set, in the limit of vanishing elastic
constant. We use this convergence result to obtain an expansion
for the scalar order parameter $s^L$ in terms of $L$, everywhere
away from the isotropic set. As mentioned in
Section~\ref{sec:prelim}, a limiting harmonic map is an energy
minimizer within the Oseen-Frank theory for uniaxial liquid
crystals with constant order parameter. These convergence results
suggest that Oseen-Frank theory and Landau-de Gennes theory give
qualitatively similar information away from topological defects
and the Landau-de Gennes theory can potentially give new
information near topological defects.

As mentioned in Section~\ref{sec:intro}, some of our results will
also apply to uniaxial solutions with bounded energy in the limit
$L\to 0^+$. In particular, the strong convergence result in
Section~\ref{sec:results1} will hold for any sequence of solutions of (\ref{eq:11}) and (\ref{eq:decoupled}) whose energy is bounded from above by the energy of the limiting harmonic map $\Q^0$, in the limit $L\to 0^+$. The interior and boundary monotonicty lemmas in Section~\ref{sec:results1} hold for all solutions of (\ref{eq:11}) (and hence, of (\ref{eq:decoupled}) which is a special case of (\ref{eq:11})). In particular, (\ref{eq:scalar}) will be valid for all sequences of uniaxial solutions, $\left\{\Q^{L_k}\right\}$ of (\ref{eq:decoupled}), whose energies are bounded above by the energy of a limiting harmonic map $\Q^0$, in the limit $L_k\to 0^+$. We will still have the $C^{1,\alpha}$-convergence of $\left\{\Q^{L_k}\right\}$ to $\Q^0$ as $L_k\to 0$, everywhere away from the singular set of the limiting harmonic map.

Such sequences of uniaxial solutions with bounded energy do exist, such as \emph{radial-hedgehog solutions} on a unit ball with strong radial anchoring conditions \cite{maj2}. Radial-hedgehog solutions are uniaxial, spherically-symmetric solutions of (\ref{eq:decoupled}). They are analogous to degree $+1$ vortices in Ginzburg-Landau theory and in \cite{maj2}, we use Ginzburg-Landau techniques to study these radial-hedgehog solutions, their defect cores and stability properties in the limit $L\to 0^+$.

\textit{Extensions to biaxial case:}
Some of the arguments in this paper can be extended to general biaxial minimizers of the form (\ref{eq:1}). As an example, let $\left\{\Q^{L_k}\right\}$ be a sequence of minimizers (biaxial or uniaxial) of the Landau-de Gennes energy functional $\Ical$, in the admissible space $\Acal_\Q$. Then $\left\{\Q^{L_k}\right\}$ converges strongly to a limiting harmonic map $\Q^0$ (as in (\ref{eq:Q0})) in $W^{1,2}(\Omega,\bar{S})$ (up to a subsequence) \cite{amaz}, for $L_k \to 0^+$ as $k\to\infty$.
Using the interior and boundary monotonicity lemmas (\ref{eq:intmon}) and (\ref{eq:bnd}), we can show that
\begin{equation}
\label{eq:biaxialnew1}
 f_B\left( \Q^{L_k}\right) \to 0 \end{equation}
uniformly everywhere away from the singular set, $S_0$, of the limiting harmonic map or equivalently
\begin{equation}
\label{eq:biaxialnew2} s\to s_+, \quad r\to 0^+
\end{equation} uniformly away from $S_0$, as $k\to \infty$ \cite{amaz}.

In what follows, we derive the analogue of Lemma~\ref{lem:3} in the biaxial case.
\begin{lemma}
\label{lem:biaxial}
Let $\Q^L$ be a global minimizer of $\Ical$, for $L$ sufficiently small \begin{footnote}{we are not making any assumptions about $\Q$; a general global minimizer for $\Ical$ in the admissible space $\Acal$ exists from the direct methods in the calculus of variations}\end{footnote}. Let $$A^L = \frac{1}{2}\Q^L_{ij,k}\Q^L_{ij,k}.$$
Then on $\Omega\setminus B_\delta(S_0)$, we have the following inequality
\begin{equation}
\label{eq:biaxial}
-\Delta A^L + \left|D^2 \Q^L \right|^2 \leq
\frac{1}{\alpha^2}\left|D^2 \Q^L \right|^2  + \alpha^4
\frac{A{^L}^2}{|\Q^L|^2}
\end{equation}
where $B_\delta(S_0)$ is a small $\delta$-neighbourhood of $S_0$ and $\alpha>1$ is a positive constant independent of $L$.
\end{lemma}

\textit{Proof:} We start with the relation (\ref{eq:rs5})
$$\Delta A^L = |D^2 \Q^L|^2 + \Q^L_{ij,k}\Q^L_{ij,ppk}$$
and drop the superscript $L$ for brevity.

We need to estimate $\left| \Q_{ij,k}\Q_{ij,ppk} \right|$ in terms of $\left|\Delta \Q\right| |\grad \Q|^2/|\Q|$.
Straightforward but tedious calculations show that
\begin{eqnarray}
&& L^2 \left| \Q_{ij,k}\Q_{ij,ppk} \right|^2 = a^4 |\grad \Q|^4 + c^4 \left(\textrm{tr}\Q^2\right)^2 |\grad \Q|^4 + 4b^4\left(\Q_{ip}\Q_{pj,q}\Q_{ij,q}\right)^2 + \nonumber \\ && + 4a^2 b^2 \Q_{ip}\Q_{pj,q}\Q_{ij,q} |\grad \Q|^2 - 4b^2 c^2 |\Q|^2 |\grad \Q|^2 \Q_{ip}\Q_{pj,q}\Q_{ij,q} - 2a^2 c^2 |\Q|^2 |\grad \Q|^4 + \nonumber \\ &&
 + 4 c^4\left(\Q\cdot\grad \Q\right)^4 + 4c^4\left(\Q \cdot\grad \Q\right)^2 |\Q|^2 |\grad \Q|^2
- 4a^2 c^2 |\grad \Q|^2\left(\Q \cdot\grad \Q\right)^2 - 8b^2 c^2 \left(\Q \cdot\grad \Q\right)^2\Q_{ip}\Q_{pj,q}\Q_{ij,q} \leq \nonumber \\ && \leq C (a^2,b^2,c^2) |\grad \Q|^4  \label{eq:biaxial2}
\end{eqnarray}
where we have used the Euler-Lagrange equations (\ref{eq:11}) to compute the right-hand side of (\ref{eq:biaxial2}) and the uniform convergence of the bulk energy density to its minimum value away from the singular set of the limiting harmonic map. It can be shown that the right-hand side of (\ref{eq:biaxial2}) vanishes for $\Q \in \Q_{\min}$, where $\Q_{\min}$ has been defined in (\ref{eq:7}). The details of these calculations are omitted here for brevity.

Secondly,
\begin{eqnarray}
&& L^2 \frac{|\grad \Q|^4}{|\Q|^2}\left| \Delta \Q \right|^2 = a^4 |\grad \Q|^4 + 2a^2 b^2 |\grad \Q|^4 \frac{\textrm{tr}\Q^3}{|\Q|^2} - 2a^2 c^2 |\Q|^2 |\grad \Q|^4 \nonumber \\ && - 2b^2 c^2 \textrm{tr}\Q^3 |\grad \Q|^4 + c^4 |\Q|^4 |\grad \Q|^4 + 2b^4 \frac{s^4 + r^4 + 3s^2r^2 - 2s^3r - 2sr^3}{27 |\Q|^2}|\grad \Q|^4
\geq D(a^2,b^2,c^2) |\grad \Q|^4\nonumber \\ && \label{eq:biaxial3}
\end{eqnarray}
where $$D(a^2,b^2,c^2) = 0$$ if and only if $\Q \in \Q_{\min}$.

Combining (\ref{eq:biaxial2}) and (\ref{eq:biaxial3}), we get that
\begin{equation}
\label{eq:biaxial4}
 \left| \Q_{ij,k}\Q_{ij,ppk} \right| \leq D^{'}(a^2,b^2,c^2)\frac{|\grad \Q|^2}{|\Q|}\left| \Delta \Q \right|
\end{equation}
where $D^{'}$ is a positive constant independent of $L$. Substituting (\ref{eq:biaxial4}) into (\ref{eq:rs5}) and repeating the same steps as in Lemma~\ref{lem:3}, (\ref{eq:biaxial}) follows. The proof of Lemma~\ref{lem:biaxial} is then complete.  $\Box$

\textit{Corollary: Let $\Q^L$ be a global minimizer of $\Ical$ (biaxial or uniaxial), in the admissible space $\Acal_\Q$, for $L$ sufficiently small. Then we have the following interior estimates, away from the singular set, $S_0$, of the limiting harmonic map $\Q^0$ in (\ref{eq:Q0}) :-
\begin{eqnarray}
&& \frac{1}{2}|\grad\Q^L|^2 + \frac{f_B(\Q^L)}{L} \leq H(a^2,b^2,c^2,\Omega) \label{eq:modQa} \\ &&
\label{eq:modQ}
|\Q^0| - |\Q^L| \leq C(a^2,b^2,c^2) L \quad on~K\subset \Omega\setminus B_\delta(S_0).
\end{eqnarray} In particular, the largest positive eigenvalue, $\la_1^L$, of $\Q^L$, satisfies the following inequality on the interior compact subset $K\subset \Omega\setminus B_\delta(S_0)$
\begin{equation}
\label{eq:modQ1}
\frac{2s_+}{3} - \la_1^L \leq D(a^2,b^2,c^2) L
\end{equation}
where $s_+$ has been defined in (\ref{eq:8}) and the positive constants $H$,$C$ and $D$ are independent of $L$.}

\textit{Proof:} The inequality (\ref{eq:modQa}) is a mere repetition of (\ref{eq:uniform3}); see \cite{amaz} for a proof.

Consider the function
$$\left| \Q^L \right| = \left(\Q^L_{pq} \Q^L_{pq}\right)^{1/2} \quad p,q=1,2,3.$$
Then a direct computation shows that $|\Q^L|$ satisfies the following partial differential equation
\begin{equation}
\label{eq:modQ2}
\Delta |\Q| = \frac{|\grad \Q|^2}{|\Q|} - \frac{\left(\Q \cdot \grad \Q\right)^2}{|\Q|^3} + \frac{\Q_{rs}\Delta\Q_{rs}}{|\Q|} \quad r,s=1\ldots 3
\end{equation}
where we have dropped the superscript $L$ for brevity.

On the interior compact subset $K\subset\Omega\setminus B_\delta(S_0)$, we have the following inequalities
\begin{eqnarray}
&& \frac{2}{3}s_+^2 - \eps_1 \leq |\Q|^2 \leq \frac{2}{3}s_+^2 \nonumber \\
&& |\grad \Q|^2 \leq C_1(a^2,b^2,c^2)
\label{eq:modQ3}
\end{eqnarray}
where we have used (\ref{eq:max}), (\ref{eq:biaxialnew1}) and (\ref{eq:modQa}).
Therefore, the first two terms on the right-hand side of (\ref{eq:modQ2}) can be bounded independently of $L$.
We use the Euler-Lagrange equations (\ref{eq:11}) to compute the third term on the right-hand side of (\ref{eq:modQ2}) i.e.
$$ \frac{\Q_{rs}\Delta\Q_{rs}}{|\Q|} =  \frac{1}{|\Q| L}\left\{-a^2 |\Q|^2 - b^2\textrm{tr}\Q^3 + c^2 |\Q|^4\right\} =  \frac{|\Q|}{ L}\left\{c^2 |\Q|^2 - \frac{b^2 |\Q|}{\sqrt{6}} - a^2\right\} + \frac{b^2 |\Q|^2}{\sqrt{6} L}\left( 1 -  \sqrt{6}\frac{\textrm{tr}\Q^3}{|\Q|^3}\right).$$

We recall from \cite{amaz} that
$$\beta^2(\Q) = 1 - 6 \left(\frac{\textrm{tr}\Q^3}{|\Q|^3}\right)^2 \in [0,1]$$
is the biaxiality parameter and as a direct consequence of (\ref{eq:modQa}), we have 
$$ \beta^2(\Q) = 1 - 6 \left(\frac{\textrm{tr}\Q^3}{|\Q|^3}\right)^2 \leq C_2(a^2,b^2,c^2) L $$
on the compact interior subset $K\subset \Omega\setminus B_\delta(S_0)$, for a positive constant $C_2$ independent of $L$. Further, we have the following sequence of inequalities on $K\subset \Omega\setminus B_\delta(S_0)$
$$ C_3(a^2,b^2,c^2) \left(|\Q| - |\Q^0|\right) \leq \left\{c^2 |\Q|^2 - \frac{b^2 |\Q|^2}{\sqrt{6}} - a^2\right\}
\leq C_4(a^2,b^2,c^2) \left(|\Q| - |\Q^0|\right)$$
for positive constants $C_3,C_4$ independent of $L$ (see (\ref{eq:biaxialnew1}) and (\ref{eq:rs1})).

From the preceding remarks, we deduce that
\begin{equation}
\label{eq:modQ4}
\Delta |\Q(\x)| \leq \alpha(a^2,b^2,c^2) +  C_4(a^2,b^2,c^2) \frac{\left(|\Q(\x)| - |\Q^0(\x)|\right)}{L},\quad \x\in \Omega\setminus B_\delta(S_0)
\end{equation}
where $\alpha$ is a positive constant independent of $L$.
Define the function $$\psi = \frac{\left(|\Q^0| - |\Q|\right)}{L}.$$ Then using (\ref{eq:modQ4}), we see that $\psi$ satisfies the following inequality on $K\subset \Omega\setminus B_\delta(S_0)$
\begin{equation}
\label{eq:modQ5}
-L \Delta \psi + \beta(a^2,b^2,c^2) \psi \leq \alpha^{'}(a^2,b^2,c^2).
\end{equation} Finally, we apply the maximum principle and Lemma~$2$ in \cite{bbh} to deduce that
\begin{equation}
\left| \psi (\x)\right| \leq \gamma(a^2,b^2,c^2) \quad \x\in \Omega\setminus B_\delta(S_0),
\end{equation} for a positive constant $\gamma$ independent of $L$ and (\ref{eq:modQ}) follows. The inequality (\ref{eq:modQ}) improves upon a previous estimate in \cite{amaz} where an analysis of the bulk energy density $f_B$, coupled with (\ref{eq:biaxialnew1}), shows that $|\Q^0| - |\Q| \leq A(a^2,b^2,c^2)\sqrt{L}$, for a positive constant $A(a^2,b^2,c^2)$ independent of $L$.

For (\ref{eq:modQ1}), we use the following alternative representation formula to (\ref{eq:1})
\begin{equation}
\label{eq:modQ6}
\Q^L = S_L\left(\n\otimes \n - \frac{1}{3}\mathbf{I}\right) + R_L\left(\m\otimes \m - \p\otimes\p \right)
\end{equation}
where $\n, \m$ and $\p$ are the orthonormal eigenvectors and
\begin{equation} \label{eq:modQ7}
0\leq s_+ - S_L \leq C_6\sqrt{L}; \quad R_L^2\leq C_5 L
\end{equation}
on the interior compact subset $K\subset \Omega\setminus B_\delta(S_0)$, for positive constants $C_6,C_5$ independent of $L$ (see (\ref{eq:biaxialnew1}) and Proposition~$7$ in \cite{amaz}). We note that $$ |\Q|^2 = \frac{2}{3}S_L^2 + 2R_L^2$$ and hence the inequality (\ref{eq:modQ}) necessarily implies that
$$ s_+ - S_L\left( 1 + \frac{3R_L^2}{S_L^2} \right)^{1/2} \leq C_7 L, $$
for a positive constant $C_7$ independent of $L$. This combined with (\ref{eq:modQ7}) i.e.  $R_L^2\leq C_5 L$ yields the improved estimate
\begin{equation}
\label{eq:modQ8}
0\leq s_+ - S_L(\x) \leq C_8 L \quad \x\in \Omega\setminus B_\delta(S_0)
\end{equation}  where $C_8>0$ is independent of $L$. Finally, it suffices to note from (\ref{eq:modQ6}) that the largest positive eigenvalue of $\Q^L$ is given by $$\la_1^L = \frac{2}{3}S_L$$ and (\ref{eq:modQ1}) directly follows from (\ref{eq:modQ8}). $\Box$

One might expect that the techniques in this paper can be generalized to obtain the analogue of Proposition~\ref{prop:convergence} and Proposition~\ref{prop:4} in the biaxial case. However, this can be accomplished only if we have a better understanding of the full Euler-Lagrange equations (\ref{eq:11}). One strategy is to decompose the system (\ref{eq:11}) as follows -
\begin{eqnarray}
\label{eq:modQ9}
&& L\Delta
\Q_{ij}=-a^2\Q_{ij}-b^2\left(\Q_{ik}\Q_{kj}-\frac{\delta_{ij}}{3}\textrm{tr}(\Q^2)\right)
+c^2\Q_{ij}\textrm{tr}(\Q^2) =  \nonumber \\ && = \left( -a^2 - b^2 \frac{|\Q|}{\sqrt{6}} + c^2 |\Q|^2\right)\Q_{ij} + b^2\left( \frac{|\Q|}{\sqrt{6}} \Q_{ij} - \Q_{ik}\Q_{kj} + \frac{1}{3}|\Q|^2\delta_{ij}\right)
\end{eqnarray}
where we can think of the first term as being a \emph{uniaxial} component and the second term as being a \emph{biaxial} component. We need to understand the coupling between the \emph{uniaxial} and the \emph{biaxial} components and to establish quantitative estimates on the magnitude of the biaxial component, in order to derive rigorous results for the structure of global Landau-de Gennes energy minimizers and their relation to the limiting harmonic map $\Q^0$ in (\ref{eq:Q0}). Other future directions are to characterize defects in Landau-de Gennes global minimizers (uniaxial versus biaxial cases), to study qualitative properties of Landau-de Gennes minimizers for different choices of the boundary conditions i.e. when $\Q_b \notin \Q_{\min}$ where $\Q_{\min}$ has been defined in (\ref{eq:7}) and to study Landau-de Gennes minimizers in different temperature regimes. We plan to report on these problems in future work.

\section{Acknowledgments}
A.~Majumdar is  supported by Award No. KUK-C1-013-04 , made by
King Abdullah University of Science and Technology (KAUST), to the
Oxford Centre for Collaborative Applied Mathematics. The author thanks Cameron Hall, Luc Ngyuen and Jonathan Robbins for helpful conversations. The author thanks Arghir Zarnescu for suggesting improvements on an earlier version of the paper.


\end{document}